\documentclass[fleqn,3p,number,11pt]{elsarticle}
\makeatletter
\def\ps@pprintTitle{%
	\let\@oddhead\@empty
	\let\@evenhead\@empty
	\def\@oddfoot{\centerline{\thepage}}%
	\let\@evenfoot\@oddfoot}
\makeatother

\usepackage{newtxtext}
\usepackage[vvarbb]{newtxmath}
\usepackage{esint}
\DeclareMathAlphabet{\mathcal}{OMS}{cmsy}{m}{n}
\usepackage{stmaryrd}
\usepackage{bm}
\usepackage{graphicx}
\usepackage[dvipsnames]{xcolor}
\usepackage{tikz}
\usepackage{pgfplots}
\usepackage{caption}
\usepackage[labelfont={bf,footnotesize},textfont=footnotesize]{subcaption}
\usepackage{adjustbox}
\usepackage[colorlinks=true]{hyperref}
\usepackage[capitalise]{cleveref}
\usepackage{braket}
\usepackage{nicefrac}
\usepackage{float}
\usepackage{multirow}
\pgfplotsset{compat=newest}
\usepackage{lineno}
\usepackage{cases}
\usepackage{booktabs}
\usepackage{dcolumn}

\makeatletter
\AtBeginDocument{\def\@citecolor{cyan}}
\AtBeginDocument{\def\@linkcolor{cyan}}
\makeatother

\newcommand{\dee}{\,\mathrm{d}}
\newcommand{\trp}{{\text{T}}}
\DeclareMathOperator{\trace}{tr}

\let\epsilon\varepsilon
\let\upepsilon\upvarepsilon

\title{{\bfseries A combined finite element-finite volume framework for phase-field fracture}}

\begin{document}
%\pagewiselinenumbers

%\author{Juan Michael Sargado \and Eirik Keilegavlen \and Inga Berre \and Jan Martin Nordbotten}

%\institute{
%	J.M. Sargado \and E. Keilegavlen \and I. Berre \and J.M. Nordbotten \at Department of Mathematics, University of Bergen \\
%	All\'egaten 41, Bergen 5007, Norway 
%	\and 
%	J.M. Sargado \at NORCE Norwegian Research Centre AS \\
%	Nyg{\aa}rdsgaten 112, Bergen 5008, Norway \\
%	\email{jusa@norceresearch.no}
%}
	

\author[uib,norce]{Juan Michael Sargado\corref{cor1}}
\ead{Juan.Sargado@uib.no}
\author[uib]{Eirik Keilegavlen}
\ead{Eirik.Keilegavlen@uib.no}
\author[uib]{Inga Berre}
\ead{Inga.Berre@uib.no}
\author[uib,pu]{Jan Martin Nordbotten}
\ead{Jan.Nordbotten@uib.no}

\cortext[cor1]{Corresponding author}
\address[uib]{Department of Mathematics, University of Bergen, All\'egaten 41, 5007 Bergen, Norway}
\address[norce]{NORCE Norwegian Research Centre AS, Nyg{\aa}rdsporten, Bergen, Norway}
\address[pu]{Princeton Environmental Institute, Princeton University, Princeton, NJ 08544, USA}

\journal{Computer Methods in Applied Mechanics and Engineering}

\begin{abstract}
\small
Numerical simulations of brittle fracture using phase-field approaches often employ a discrete approximation framework that applies the same order of interpolation for the displacement and phase-field variables. Most common is to use linear finite elements to discretize the linear momentum and phase-field equations. However the use of $P_1$ Lagrange shape functions to model the phase-field is not optimal, since the latter develops cusps for fully developed cracks that in turn occur at locations correspoding to Gauss points of the associated FE model for the mechanics. Such feature is challenging to reproduce accurately with low order elements, and consequently element sizes must be made very small relative to the phase-field regularization parameter in order to achieve convergence of results with respect to the mesh. In this paper, we combine the standard $P_1$ FE discretization of stress equilibrium with a cell-centered finite volume approximation of the phase-field evolution equation based on the two-point flux approximation that is constructed on the same simplex mesh. Compared to a pure FE formulation utilizing linear elements, the proposed framework results in looser restrictions on mesh refinement with respect to the phase-field length scale. Furthermore, initialization of the history field is straightforward and accomplished through a local procedure. The ability to employ a coarser mesh relative to the traditional implementation is shown for several numerical examples, demonstrating savings in computational cost on the order of 50 to 80 percent for the studied cases.
\end{abstract}


\maketitle

\section{Introduction}
Brittle fracture is an important failure mechanism against which engineering structures must be properly designed to ensure their safety. At the same time, the deliberate initiation and propagation of fractures are central to many subsurface applications related to energy production and waste disposal. Griffith \cite{Griffith1921} was the first to formulate a theory of brittle fracture based on thermodynamic arguments, utilizing earlier results by Inglis \cite{Inglis1913} concerning stresses around elliptical holes. Griffith's theory was later amended to include the effect of plastic zones around crack tips by Irwin, who also introduced the notion of stress intensity factors \cite{Irwin1957}. Together, these form the basis of classical linear elastic fracture mechanics (LEFM). More recently, an extension to LEFM in the form of a variational theory of fracture was formulated by Francfort and Marigo \cite{Francfort1998}, with the aim of overcoming the former's limitations with regard to crack initiation and branching through the adoption of an energy minimization framework.

Meanwhile, establishment of the finite element method in the early 1960s also fueled the development of numerical techniques for simulating fracture in solid structures. In particular, the work of Ngo and Scordelis \cite{Ngo1967} and Rashid \cite{Rashid1968} on crack formation in concrete are recognized as pioneering efforts in discrete and continuous approaches to fracture modeling. In the former, cracks are treated as discrete/sharp entities, and the introduction of new crack segments can be realized explicitly via inter-element separation algorithms \cite{Xu1994,Camacho1996}, ideally in combination with adaptive remeshing in order to remove mesh bias. As brittle material behavior is characterized by the occurrence of stress singularities at crack tips, this necessitates considerable mesh refinement in order to accurately resolve local stress fields. A classical approach of dealing with this problem is to make use of quarter-point elements \cite{Barsoum1977}.  Alternatively, an implicit treatment of discrete cracks can be made through enrichment of the approximation space by suitable functions that directly model displacement discontinuities and stress singularities, as in the extended finite element method \cite{Belytschko1999,Moes1999}. On the other hand in smeared crack approaches, no effort is made to track the fracture surfaces. Instead, material stiffness is progressively reduced in order to approximate the overall mechanical response arising from the presence of a crack. Notable developments in this category include the fictitious crack model \cite{Hillerborg1976}, rotating smeared cracks \cite{Gupta1984}, the crack band model \cite{Bazant1983}, integral-type damage theory \cite{Pijaudier1987}, implicit gradient-enhanced damage models \cite{Peerlings1996,Peerlings1998}, and variational phase-field models \cite{Bourdin2000,Bourdin2008}. The last in particular have been gaining massive popularity in recent years due to their clear connection to Griffith's theory. Phase-field models fall under the continuous approach to fracture and work by representing cracks as diffuse entities through an auxiliary scalar variable known as the phase-field \cite{Kuhn2010}, with the amount of diffusion being controlled by a regularization parameter $\ell$. While the initial intent of \cite{Bourdin2000} was for the regularized energy functional to approximate a body with discrete cracks through the concept of $\Gamma$-convergence \cite{Braides2006}, lately it has come to be understood that variational phase-field models are in fact a subset of gradient-enhanced damage models, and that the phase-field regularization parameter should be related to the intrinsic length scale of the material \cite{Pham2011,Sicsic2013,Bourdin2014,Tanne2018}.

\begin{figure}
	\centering
	\begin{subfigure}[b]{0.32\textwidth}
		\includegraphics[width=\textwidth]{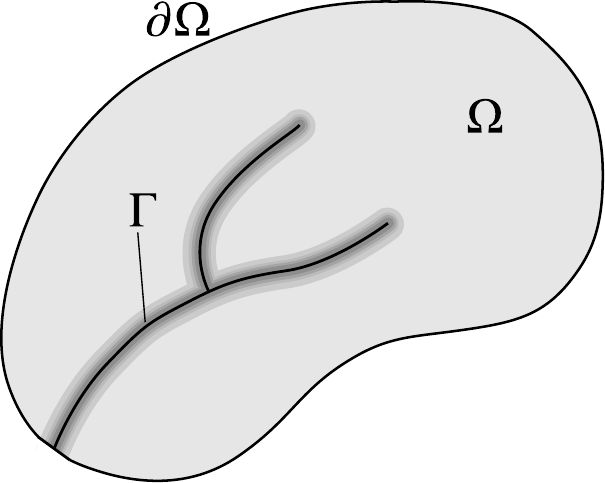}
		\caption{}
	\end{subfigure} \hspace{10pt}
	\begin{subfigure}[b]{0.45\textwidth}
		\includegraphics[width=\textwidth]{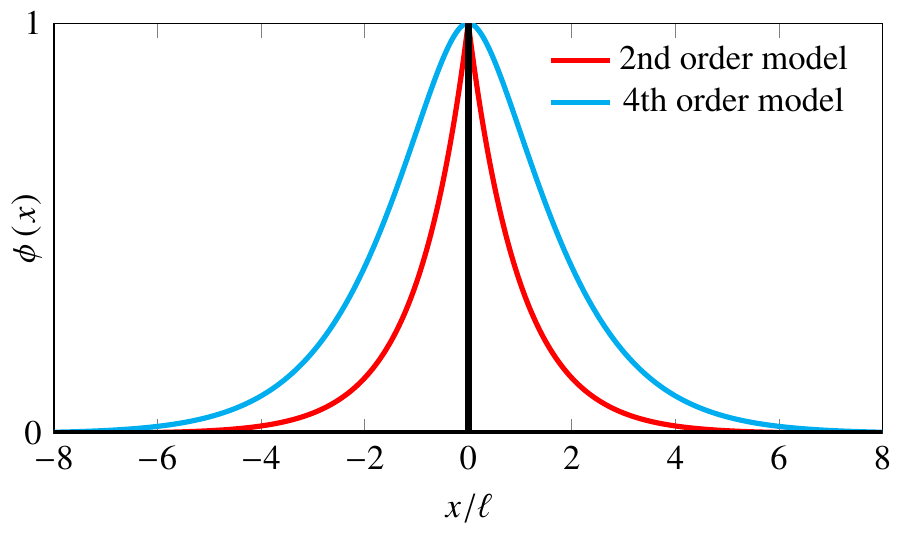}
		\caption{}
	\end{subfigure}
	\caption{Phase-field modeling of fractures. (a) Discrete cracks $\Gamma$ are assumed to be surfaces of dimension $n-1$ embedded inside an $n$-dimensional body; these are represented as diffuse entities via an auxiliary scalar field $\phi \in \left[ 0,1 \right]$. (b) Phase-field profile for a 1-dimensional body with a fully developed crack at $x = 0$. The classical 2nd order phase-field model results in $\phi$ having a cusp at the crack location, which can be rounded by using a higher order model, e.g.\ the 4th order model proposed in  \cite{Borden2014}.}
	\label{fig:phaseFieldRegurization}
\end{figure}

The main strength of the phase-field method relative to other approaches is the ease with which complex evolution of fractures can be modeled. In particular, its underlying energy minimization framework automatically handles the branching of fractures without the need to introduce external criteria. Furthermore since cracks are modeled by a single field, no additional complexity involving bookkeeping of surface intersections and the like is introduced upon the initiation and branching/coalescence of fractures. Instead, an additional partial differential equation must be solved that governs the auxiliary field evolution. This makes the approach promising for multiphysics applications where the interaction between different physical processes due to fracture formation can be recast in a diffuse setting, essentially becoming a coupling term that is additionally dependent on the phase-field. Such applications have included thermal cracking \cite{Bourdin2010,Bourdin2014}, electro-mechanical processes \cite{Abdollahi2011,Wilson2013}, fluid-driven fracture in elastic and poroelastic materials \cite{Mikelic2015_nonlin,Mikelic2015_cg,Miehe2015,Miehe2016_cmame,Yoshioka2016,Cajuhi2017}, and chemo-mechanical degradation of battery electrode particles \cite{Klinsmann2016,Miehe2016_ijnme}. On the other hand, a prevailing challenge in the use of fracture phase-field models is the computational expense associated with their solution, primarily due to the nonlinear coupling of the component PDEs and also the need to properly resolve the diffuse cracks in the mesh according to the phase-field length scale. While the crack tip stress singularities associated with discrete approaches are eliminated in the phase-field formulation \cite{Sicsic2013} (and with them the need to aggressively refine meshes according to the stress distribution), selective mesh refinement is nevertheless required as the phase-field profile for fully developed cracks contains cusps as illustrated in Figure \ref{fig:phaseFieldRegurization}.
Consequently the discretization must be suitable fine in the vicinity of cracks to resolve rapidly changing gradients. This can dramatically increase the computational size of the discrete problem, particularly if the crack path is unknown a priori and a global mesh refinement strategy is used.

Remarkably, the above nuance involving mesh size requirements (i.e.\ that they emanate primarily from the phase-field subsystem rather than mechanics) has not been exploited in prior literature. Rather, the standard practice thus far has been to employ the same order of basis functions for the linear momentum and phase-field equations. In addition, while earlier studies utilizing FEM such as \cite{Kuhn2010,Miehe2010_ijnme} favored setting the characteristic size of elements to half the magnitude of $\ell$ in conjunction with $P_1$ basis functions, this ratio has been found to be suboptimal in later studies \cite{Nguyen2015,Sargado2018}. Higher order models have also been introduced which result in looser requirements on the mesh resolutions with respect to $\ell$ due to blunting of the cusp, however such models require $C^1$-continuous solutions for the phase-field and hence are more suited for use with specialized schemes that preserve higher order continuity, for instance isogeometric basis functions \cite{Borden2014} and local maximum entropy approximants \cite{Li2015}. Another possible means of reducing the computational burden is to perform adaptive remeshing as demonstrated in \cite{Borden2012,Klinsmann2015}, with the caveat that such procedures also come with their own added computational cost. In particular, changes in the data structure in the form of addition/removal of cells and nodes can turn out to be more expensive than plain assembly and solution of non-adaptive problems that properly exploit the static nature of the sparse matrix profiles. An interesting direction recently explored in \cite{Giovanardi2017} combines discrete approaches (in the form of extended finite elements) with the phase-field method, in which the latter is restricted to local domains surrounding crack tips in order to limit the overall size of the discrete problem.

In this contribution, we propose a combined discretization scheme (in the sense of both unequal order approximation and different formulations) to solve the coupled system of equations associated with phase-field modeling of brittle fracture. This idea is based on the assumption that the need to properly resolve the phase-field along cracks results in mesh size restrictions that are more stringent than what is needed for the accurate resolution of stresses. Within a finite element context, this naturally leads to the prospect of using higher order basis functions for the auxiliary field. Alternatively, one can make use of an approximation scheme that is better equipped to handle the unique challenges of the phase-field method. In particular for this paper, we look at the prospect of adopting control volume formulation approaches such as finite volumes. Compared with weak formulation-based methods such as FEM, the use of finite volumes in solid mechanics applications has not been as extensively explored, and to our knowledge the only application thus far of the latter in phase-field fracture has been the work of Santillan et al.\ \cite{Santillan2017}, who utilized a finite volume discretization for all of the governing equations. In this work, we propose instead to combine a finite element approximation of the linear momentum equation using $P_1$ shape functions with a cell-centered finite volume scheme for the phase-field subsystem based on the classical two-point flux approximation (TPFA). The phase-field is assumed to be piecewise constant over elements, hence our scheme is effectively a $P_1$-$P_0$ formulation. However the main advantage of a finite volume scheme in this case is that it allows for discontinuous gradients within a local cell while preserving the continuity of normal derivative terms across cell faces, which is ideal for modeling the cusps that occur in the phase-field profile for fully developed cracks. Furthermore, local calculations pertaining to a cell do not involve matrix operations and hence are much cheaper to carry out than those associated with a variational formulation such as finite elements. A limitation of the TPFA is that it ignores tangential components of gradients across cell faces, hence in porous media flow applications for which the technique was originally developed, the cell faces should be aligned with the principal directions of the permeability tensor in order to produce correct results \cite{Aavatsmark1996}. On the other hand in the PDE governing fracture evolution, the relevant material parameter acting on the phase-field gradient is the critical energy release rate, $\mathcal{G}_c$. For a wide class of materials, $\mathcal{G}_c$ is specified as a scalar quantity that can alternatively be interpreted as an isotropic tensor, and our experience in this case has been that application of TPFA to unstructured simplex meshes yields acceptable results despite cell faces being generally non-orthogonal.
\section{Mathematical formulation}
Let us consider a solid body occupying a domain $\Omega$ and containing a collection of internal cracks denoted by $\Gamma$. Following \cite{Francfort1998}, we assume that its total potential energy may be expressed as a sum of bulk and surface terms. That is,
\begin{linenomath}
\begin{equation}
	\Psi = \int_{\Omega\setminus\Gamma} \psi_0 \left( \bm{\upepsilon} \right) \dee\Omega + \int_\Gamma \mathcal{G}_c \dee\Gamma,
	\label{eq:sharpCrackFunctional}
\end{equation}
\end{linenomath}
wherein $\psi_0 \left( \bm{\upepsilon} \right) = \frac{1}{2} \lambda \left( \trace \bm{\upepsilon} \right)^2 + \mu \bm{\upepsilon} : \bm{\upepsilon}$ is the Helmholtz free energy for the bulk material, and $\mathcal{G}_c$ is the critical energy release rate from Griffith's theory of brittle fracture \cite{Griffith1921}. As the second term represents the energy dissipation associated with crack formation, it is further subject to the irreversibility constraint
\begin{linenomath}
\begin{equation}
	\Gamma \left( t + \Delta t \right) \supseteq \Gamma \left( t \right)
\end{equation}
\end{linenomath}
for $\Delta t > 0$, meaning that cracks may only grow over time but not heal. Determining the proper evolution of $\Gamma$ under arbitrary loading conditions is nontrivial, and following \cite{Bourdin2008} we instead utilize a regularized version of the above potential, expressed as
\begin{linenomath}
\begin{equation}
	\Psi_\ell = \int_\Omega \psi \left( \bm{\upepsilon}, \phi \right) \dee\Omega + \int_\Omega \mathcal{G}_c \gamma_\ell \left( \phi, \nabla\phi \right) \dee\Omega.
\end{equation}
\end{linenomath}
In the above expression, the cracks have been recast as diffuse entities, with $\gamma_\ell \left( \phi, \nabla\phi \right)$ representing the crack surface density per unit volume, or analogously the crack length per unit area for 2-dimensional problems. In particular we adopt a form derived in \cite{Miehe2010_ijnme} and given by
\begin{linenomath}
\begin{equation}
	\gamma_\ell \left( \phi, \nabla\phi \right) = \frac{1}{2\ell} \phi^2 + \frac{\ell}{2} \nabla\phi \cdot \nabla\phi .
	\label{eq:regularizedCrack}
\end{equation}
\end{linenomath}
The scalar variable $\phi \in \left[ 0,1 \right]$ is the known as the crack phase-field, and characterizes respectively for $\phi = 0$ and $\phi = 1$ the fully pristine and fully broken states. On the other hand, $\ell$ is a regularization parameter (the phase-field length scale) that controls the amount of diffusion in the crack representation.

There are several ways to regularize the bulk term in \eqref{eq:sharpCrackFunctional}. In the present study, we make use of two different regularizations: the isotropic model
\begin{linenomath}
\begin{equation}
	\psi \left( \bm{\upepsilon}, \phi \right) = g \left( \phi \right) \psi_0 \left( \bm{\upepsilon} \right)
	\label{eq:bourdinModel}
\end{equation}
\end{linenomath}
originally employed in \cite{Bourdin2000}, and an anisotropic model for approximating unilateral contact developed by Amor et al.\ \cite{Amor2009}, given by
\begin{linenomath}
\begin{equation}
	\begin{split}
		\psi \left( \bm{\upepsilon}, \phi \right) &= g \left( \phi \right) \left[ \frac{\kappa}{2} \left( \trace^+ \bm{\upepsilon} \right)^2 + \mu \bm{\upepsilon}_\text{dev} : \bm{\upepsilon}_\text{dev} \right] + \frac{\kappa}{2} \left( \trace^- \bm{\upepsilon} \right)^2 \\
		&= g \left( \phi \right) \psi_0^+ \left( \bm{\upepsilon} \right) + \psi_0^- \left( \bm{\upepsilon} \right)
	\end{split}
	\label{eq:amorModel}
\end{equation}
\end{linenomath}
in which $\bm{\upepsilon}_\text{dev}$ denotes the deviatoric component of $\bm{\upepsilon}$ (see \ref{sec:amorDetails}), and $\trace^\pm \bm{\upepsilon}$ is defined as
\begin{linenomath}
\begin{equation}
	\begin{split}
	\trace^+ \bm{\upepsilon} &= \max \left( 0, \trace \bm{\upepsilon} \right) \\
	\trace^- \bm{\upepsilon} &= \min \left( 0, \trace \bm{\upepsilon} \right).
	\end{split}
\end{equation}
\end{linenomath}
Both models above involve an energy degradation function $g \left( \phi \right)$, whose role is to annihilate material stiffness at locations where the material is broken according to the phase-field. For simplicity, we adopt the quadratic expression
\begin{linenomath}
\begin{equation}
	g \left( \phi \right) = \left( 1 - \phi \right)^2
\end{equation}
\end{linenomath}
that is most often used in the literature. Nevertheless we note that the above form is not the most optimal for use with \eqref{eq:regularizedCrack}, and in particular for accurate modeling of failure loads as well as overall linear response prior to fracture we refer the reader to our previous work on parametric degradation functions \cite{Sargado2018}.

To obtain the governing equations for the brittle fracture problem, we first define a functional that includes both the potential energy and external work terms, i.e.
\begin{linenomath}
\begin{equation}
	\Pi = \Psi_\ell - \int_\Omega \mathbf{b} \cdot \mathbf{u} \dee\Omega - \int_{\partial\Omega^N} \mathbf{t} \cdot \mathbf{u} \dee\partial\Omega
	\label{eq:TotalFunctional}
\end{equation}
\end{linenomath}
in which $\mathbf{b}$ denotes the body force, and $\mathbf{t}$ the surface traction acting on the Neumann boundary $\partial\Omega^N$. Imposing stationarity of the above functional, we have
\begin{linenomath}
\begin{equation}
	\delta\Pi = 0 = \frac{\partial\Pi}{\partial\mathbf{u}} \delta\mathbf{u} + \frac{\partial\Pi}{\partial\phi} \delta\phi
	\label{eq:stationarity}
\end{equation}
\end{linenomath}
which yields two coupled PDEs corresponding to the different terms in the right hand side above. The first term is none other than the weak form of the stress equilibrium equation:
\begin{linenomath}
\begin{equation}
	\int_\Omega \bm{\upsigma} \left( \bm{\upepsilon}, \phi \right) : \delta\bm{\upepsilon} \dee\Omega = \int_\Omega \mathbf{b} \cdot \delta\mathbf{u} \dee\Omega + \int_{\partial\Omega^N} \mathbf{t} \cdot \delta\mathbf{u} \dee\partial\Omega ,
	\label{eq:linearMomentumWeakForm}
\end{equation}
\end{linenomath}
wherein $\bm{\upsigma} \left( \bm{\upepsilon}, \phi \right) = \partial \psi \left( \bm{\upepsilon}, \phi \right) / \partial\bm{\upepsilon}$. On the other hand, expanding the second term in \eqref{eq:stationarity} yields the weak form of the phase-field evolution equation, given by
\begin{linenomath}
\begin{equation}
	\int_\Omega \mathcal{G}_c \left( \frac{1}{\ell} \phi \,\delta\phi + \ell \nabla\phi \cdot \nabla\delta\phi \right) \dee\Omega + \int_\Omega g^\prime \left( \phi \right) \psi_0^+ \left( \bm{\upepsilon} \right) \dee\Omega = 0
	\label{eq:phaseFieldEvolutionWeakForm}
\end{equation}
\end{linenomath}
As in \cite{Bourdin2000}, we impose homogeneous natural conditions on the entire external boundary, i.e.
\begin{linenomath}
\begin{equation}
	\nabla\phi \cdot \mathbf{n} = 0 \quad \text{on } \partial\Omega,
\end{equation}
\end{linenomath}
where $\mathbf{n}$ is the outward unit normal vector to $\partial\Omega$. Performing integration by parts and imposing arbitrariness of $\delta\mathbf{u}$ and $\delta\phi$, we obtain the strong form of coupled system given by the following boundary value problem:
\begin{linenomath}
\begin{numcases}{}
	\nabla \cdot \bm{\upsigma} \left( \bm{\upepsilon}, \phi \right) + \mathbf{b} = \mathbf{0} & in $\Omega$
	\label{eq:linearMomentumStrongForm} \\
	\mathbf{u} = \bar{\mathbf{u}} & on $\partial\Omega^D$ \\
	\bm{\upsigma} \left( \bm{\upepsilon}, \phi \right) \cdot \mathbf{n} = \mathbf{t} & on $\partial\Omega^N$ \\
	\mathcal{G}_c \ell \nabla^2 \phi - \frac{\mathcal{G}_c}{\ell} \phi = g^\prime \left( \phi \right) \psi_0^+ \left( \bm{\upepsilon},\phi \right) & in $\Omega$
	\label{eq:phaseFieldStrongForm} \\
	\nabla\phi \cdot \mathbf{n} = 0 \quad & on $\partial\Omega$
\end{numcases}
\end{linenomath}
It then remains to enforce the irreversibility of crack growth. Here we adopt the approach of Miehe et al.\ \cite{Miehe2010_cmame} and replace $\psi_0^+$ with a history field defined as
\begin{linenomath}
\begin{equation}
	\mathcal{H} \left( \mathbf{x}, t \right) = \left\{ \begin{array}{ll}
	\max\limits_{s \in \left[ 0,t \right]} \psi_0^+ \left( \bm{\upepsilon} \left( \mathbf{x}, s \right) \right), & \phi > \phi_c \\[10pt]
	\psi_0^+ \left( \bm{\upepsilon} \left( \mathbf{x}, t \right) \right) & \text{otherwise}.
	\end{array} \right.
	\label{eq:historyField}
\end{equation}
\end{linenomath}
in which $\phi_c$ is some chosen irreversibility threshold.
\section{Computational framework}
\subsection{Discrete approximation}
\begin{figure}
	\centering
	\begin{subfigure}[b]{0.35\textwidth}
		\includegraphics[width=\textwidth]{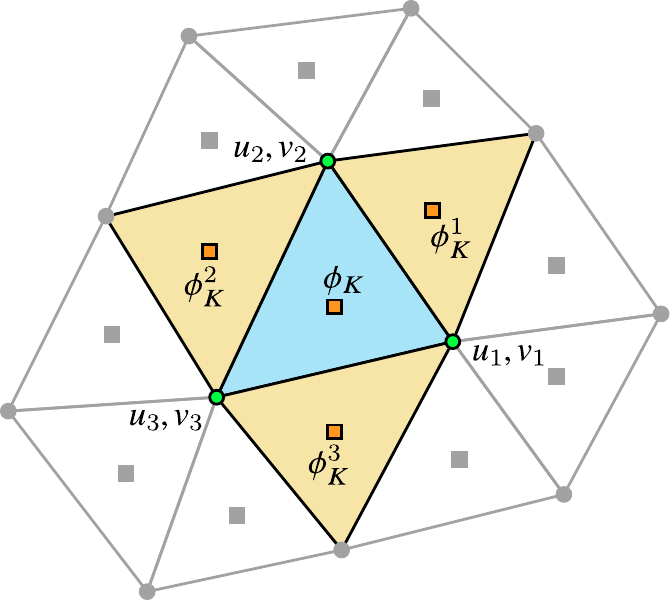}
		\caption{}
		\label{fig:stencil}
	\end{subfigure} \hspace{50pt}
	\begin{subfigure}[b]{0.25\textwidth}
		\includegraphics[width=\textwidth]{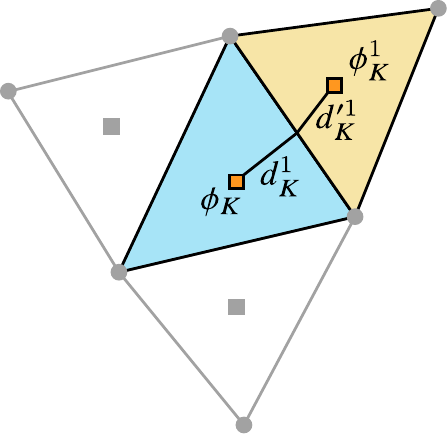} \vspace{5pt}
		\caption{}
		\label{fig:tpfa}
	\end{subfigure}
	\caption{Discrete formulation of the phase-field brittle fracture problem, combining $P_1$ finite elements for the linear momentum equation and cell-centered finite volumes for the phase-field evolution equation. (a) Computational stencil for $\Omega_K$, and (b) distances from respective cell centers to the face midpoint that are used to calculate the transmissibility coefficient $M_K^1$ between $\Omega_K$ and $\Omega_K^1$.}
\end{figure}

To perform a numerical solution of the coupled system, we combine a finite element formulation of the linear momentum equation using $P_1$ basis functions with a cell-centered finite volume discretization of the phase-field equation based on the two-point flux approximation. Figure \ref{fig:stencil} shows the resulting computational stencil for an interior triangular cell $\Omega_K$, which involves displacement degrees of freedom at the cell vertices and phase-field DOFs located at cell centers of $\Omega_K$ and its immediate surrounding cells, denoted by $\Omega_K^i$.
To construct the local matrices pertaining to a given cell, we begin by expressing the displacement and strain fields in Voigt form as vector quantities, i.e.\ $\mathbf{u} = \left\{ u_x, u_y \right\}^\trp$ and $\bm{\upepsilon} = \left\{ \epsilon_{xx}, \epsilon_{yy}, \epsilon_{zz}, \gamma_{xy} \right\}^\trp$. These are then interpolated over a given element as
\begin{linenomath}
\begin{equation}
	\begin{split}
	&\mathbf{u} = \sum\limits_{I=1}^3 \mathbf{N}_I \, \mathbf{u}_I \qquad \delta\mathbf{u} = \sum\limits_{I=1}^3 \mathbf{N}_I \, \delta\mathbf{u}_I \\
	&\bm{\upepsilon} = \sum\limits_{I=1}^3 \mathbf{B}_I \, \mathbf{u}_I \qquad \delta\bm{\upepsilon} = \sum\limits_{I=1}^3 \mathbf{B}_I \, \delta\mathbf{u}_I ,
	\end{split}
\end{equation}
\end{linenomath}
in which the matrices $\mathbf{N}_I$ and $\mathbf{B}_I$ are given by
\begin{linenomath}
\begin{equation}
	\begin{split}
	\mathbf{N}_I = \left[ \begin{array}{cc} N_I & 0 \\[3pt] 0 & N_I \end{array} \right] \qquad
	\mathbf{B}_I = \left[ \begin{array}{cc} N_{I,x} & 0 \\[3pt] 0 & N_{I,y} \\[3pt] 0 & 0 \\[3pt] N_{I,y} & N_{I,x} \end{array} \right]
	\end{split}
\end{equation}
\end{linenomath}
corresponding to plane strain conditions. $N_I$ and $\mathbf{u}_I$ are respectively the standard Lagrange ($P_1$) shape function and displacement vector associated with node $I$, and $N_{I,x}$ and $N_{I,y}$ are the derivatives of $N_I$ with respect to $x$ and $y$. Let us consider a local element $\Omega_K$ having area $A_K$. The damaged-reduced stress $\bm{\sigma} \left( \bm{\upepsilon}, \phi \right)$ is assumed to be constant over $\Omega_K$, hence the local residual for the discretized linear momentum equation associated with node $I$ can be integrated using one Gauss point located at the center of $\Omega_K$. This yields
\begin{linenomath}
\begin{equation}
		\mathbf{r}_I^\mathbf{u} = A_K \mathbf{B}_I^\trp \bm{\upsigma} \left( \bm{\upepsilon}, \phi \right) - A_K \mathbf{N}_I^\trp \mathbf{b} - \int_{\left( \partial\Omega_\mathbf{u}^N \right)_K} \mathbf{N}_I^\trp \mathbf{t} \dee S .
		\label{eq:mechanicsResidual}
\end{equation}
\end{linenomath}
Meanwhile, we construct the discrete approximation of the residual pertaining ot the phase-field equation by first integrating \eqref{eq:phaseFieldStrongForm} over $\Omega_K$ and applying integration by parts to obtain the corresponding control volume formulation:
\begin{linenomath}
\begin{align}
	0 &= \int_{\Omega_K} \left( g^\prime \left( \phi \right) \mathcal{H} + \frac{\mathcal{G}_c}{\ell} \phi - \mathcal{G}_c \ell \nabla^2 \phi \right) \dee\Omega \nonumber \\
	&= \int_{\Omega_K} \left[ g^\prime \left( \phi \right) \mathcal{H} + \frac{\mathcal{G}_c}{\ell} \phi - \nabla \cdot \left( \mathcal{G}_c \ell \nabla\phi \right) \right] \dee\Omega \nonumber \\
	&= \int_{\Omega_K} \left[ g^\prime \left( \phi \right) \mathcal{H} + \frac{\mathcal{G}_c}{\ell} \phi_K \right] \dee\Omega - \int_{\partial\Omega_K} \mathcal{G}_c \ell \nabla\phi \cdot \mathbf{n} \dee\partial\Omega_K \label{eq:controlVolumeForm_phaseField}
\end{align}
\end{linenomath}
Assuming $\phi$ to be piecewise constant over $\Omega_k$, the residual for the phase-field equation can be written as
\begin{linenomath}
\begin{equation}
	r_K^\phi = A_K \left[ g^\prime \left( \phi_K \right) \mathcal{H}_K + \frac{{\mathcal{G}_c}_K}{\ell_K} \phi_K \right] + \sum\limits_{i=1}^3 M_K^i \left( \phi_K - \phi_K^i \right)
	\label{eq:phaseFieldResidual}
\end{equation}
\end{linenomath}
wherein $A_K$ is the cell area, $\phi_K$ the value of phase-field inside $\Omega_K$, and ${\mathcal{G}_c}_K$ and $\ell_K$ are respectively the values of the critical energy release rate and phase-field regularization parameter within $\Omega_K$. Meanwhile $\phi_K^i$ denotes the phase-field value at the adjacent cell sharing edge $i$ with cell $K$, and $M_K^i$ is the transmissibility coefficient associated with the edge $\partial\Omega_K^i$, defined as
\begin{linenomath}
\begin{equation}
	M_K^i = \frac{L_i}{\cfrac{d_K^i}{{\mathcal{G}_c}_K \ell_K} + \cfrac{d^{\prime i}_k}{{\mathcal{G}_c}_K^i \ell_K^i}}
\end{equation}
\end{linenomath}
in which $L_i$ is the length of edge $i$, and $d_K^i$ and $d^{\prime i}_K$ are the respective distances from the centers of $\Omega_K$ and $\Omega_K^i$ to the midpoint of $\partial\Omega_K^i$ as illustrated in Figure \ref{fig:tpfa}. We note that \eqref{eq:phaseFieldResidual} is valid only when $\Omega_K$ is an interior cell. For cells situated at the boundary, one must account for prescribed boundary conditions acting on cell edges.

\subsection{Nonlinear solution}
\begin{figure}
	\centering
	\includegraphics[width=0.6\textwidth]{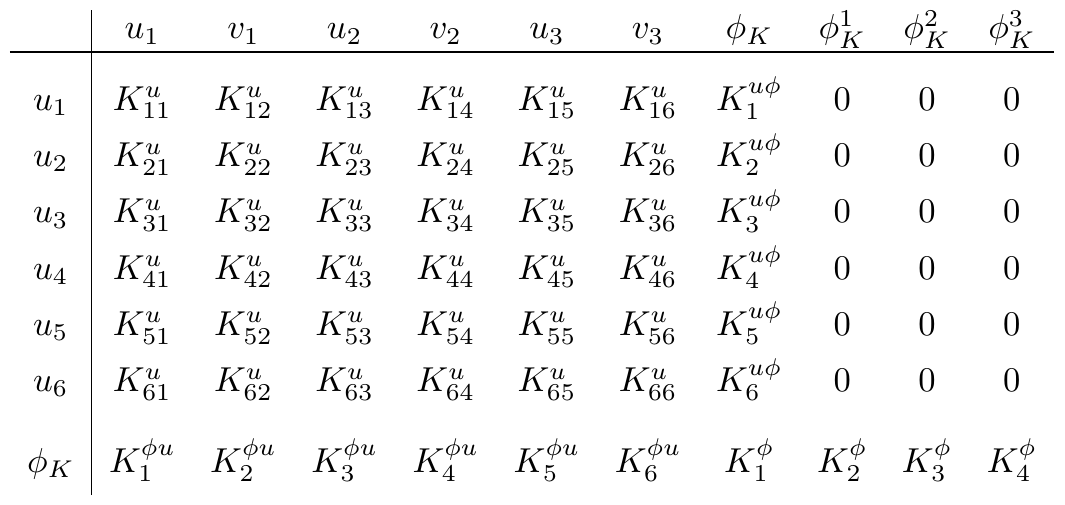}
	\caption{Local coefficient matrix associated with the proposed FE-FV scheme.}
	\label{fig:localMatrix}
\end{figure}

The system unknowns $\left\{ \mathbf{u}_I \right\}$ and $\left\{ \phi_K \right\}$ are obtained by enforcing the condition
\begin{linenomath}
\begin{numcases}{}
	\mathbf{r}_I^\mathbf{u} = \mathbf{0} \\
	r_K^\phi = 0
	\label{eq:nonlinearSystem}
\end{numcases}
\end{linenomath}
for each node $I$ and cell $K$ of the discretized domain. Due to the forms of $\bm{\upsigma} \left( \bm{\upepsilon},\phi \right)$ and $g \left( \phi \right)$, the global system is nonlinear and must be solved iteratively. Linearization of the local system yields
\begin{linenomath}
\begin{equation}
	\left[ \kern-2pt \begin{array}{ll} \mathbf{K}^u & \mathbf{K}^{u\phi} \\
		\mathbf{K}^{\phi u} & \mathbf{K}^\phi \end{array} \kern-2pt \right]
	\left\{ \kern-2pt \begin{array}{c} \dot{\mathbf{u}} \\ \dot{\bm{\upphi}} \end{array} \kern-2pt \right\} =
	\left\{ \kern-2pt \begin{array}{c} \mathbf{r}^u \\ r^\phi \end{array} \kern-2pt \right\}
	\label{eq:localSystem}
\end{equation}
\end{linenomath}
in terms of the corrections $\dot{\mathbf{u}}$ and $\dot{\bm{\upphi}}$, wherein
\begin{linenomath}
\begin{align}
	\mathbf{K}_{IJ}^u &= A_K \mathbf{B}_I^\trp \mathbb{C} \left( \bm{\upepsilon}, \phi \right) \mathbf{B}_J \\
	\mathbf{K}_I^{u\phi} &= A_K \mathbf{B}_I^\trp \frac{\partial\bm{\upsigma} \left( \bm{\upepsilon}, \phi \right)}{\partial\phi} \\
	\mathbf{K}_J^{\phi u} &= {\mathbf{K}_I^{u\phi}}^\trp \\
	\mathbf{K}^\phi &= \left[ \begin{array}{cccc} A_K \left( g^{\prime\prime} \left( \phi_k \right) \mathcal{H}_K + \dfrac{{\mathcal{G}_c}_K}{\ell_K} \right) + \sum\limits_{i=1}^3 M_K^i, & -M_K^1, & -M_K^2, & -M_K^3 \end{array} \right]
\end{align}
\end{linenomath}
Note that while the resulting global set of equations involves a symmetric coefficient matrix, the local system given by \eqref{eq:localSystem} is not square as illustrated in Figure \ref{fig:localMatrix}.
Due to the non-convexity of the functional \eqref{eq:TotalFunctional} with respect to the pair of arguments $\left( \mathbf{u}, \phi \right)$, monolithic solutions of \eqref{eq:nonlinearSystem} based on a naive application of Newton's method often fail to achieve convergence at time steps in which cracks are evolving. In the present work, we employ the alternate minimization algorithm \cite{Bourdin2008} where in each iteration we first solve for the displacement unknowns with the phase-field fixed. Said displacements are then updated and fixed, after which the proper values of the history field are computed and the phase-field subsystem solved and updated. This is carried out repeatedly until the corrections to $\mathbf{u}$ and $\phi$ as well as the residuals for each subsystem are within specified tolerances.

\subsection{Initialization of history field} \label{sec:histFieldInitialization}
Many problems of interest feature preexisting cracks in the domain. This is especially true for subsurface applications where rock masses are often characterized by ubiquitous natural fractures. In such cases, the history field of \eqref{eq:historyField} must be initialized to account for preexisting cracks. In a variational formulation-based numerical method such as finite elements, this procedure is nontrivial as the number of local equations to be satisfied is equal to the amount of phase-field DOFs per element while the number of local $\mathcal{H}$ unknowns is essentially equal to the number of integration points. Consequently, it is common to use ad hoc procedures to accomplish the initialization as outlined for example in \cite{Borden2012,Klinsmann2015}. These essentially involve performing global searches and nearest-distance calculations, and furthermore do not elaborate on how to treat locations in the immediate vicinity of multiple intersecting cracks. On the other hand in the cell-centered FV formulation employed in the present work, each local cell has just one history field unknown and also one equation pertaining to the phase-field. Hence back-calculation of $\mathcal{H}$ is straightforward which allows for initialization of the history field to be done via a two-step process:
\begin{enumerate}[(a)]
	\setlength{\itemsep}{0pt}
	\setlength{\parskip}{0pt}
	\item Apply the constraint $\phi = 1$ on all cells containing segments of the initial cracks and solve the control volume approximation of \eqref{eq:phaseFieldStrongForm} with homogeneous right hand sides, i.e.
	\begin{linenomath}
	\begin{equation}
		A_K \frac{{\mathcal{G}_c}_K}{\ell_K} \phi_K + \sum\limits_{i=1}^3 M_K^i \left( \phi_K - \phi_K^i \right) = 0
		\label{eq:histField_step1}
	\end{equation}
	\end{linenomath}
	\item Remove constraints and back-calculate the history field as follows:
	\begin{linenomath}
	\begin{equation}
		\mathcal{H}_K = \frac{-1}{g^\prime \left( \phi_K \right)} \left[ \frac{1}{A_K} \sum\limits_{i=1}^3 M_K^i \left( \phi_K - \phi_K^i \right) + \frac{{\mathcal{G}_c}_K}{\ell_K} \phi_K \right]
		\label{eq:histField_step2}
	\end{equation}
	\end{linenomath}
\end{enumerate}
It is easy to see that the second step only applies to cells which have been previously constrained, since the right hand side in \eqref{eq:histField_step2} will be identically zero in unconstrained cells due to \eqref{eq:histField_step1}. Note however that it is not possible to impose $\phi = 1$ exactly, since $g^\prime \left( 1 \right) = 0$ resulting in division by zero during back-calculation of $\mathcal{H}$. Instead, we can set $\phi$ to be sufficiently close to 1 (e.g. $\phi = 0.99$) such that the resulting material stiffness becomes effectively negligible compared to the original undamaged response. Alternatively, one can simply skip step (b) above, and instead retain the constraints imposed in step (a) throughout the simulation.

\subsection{Implementation aspects}
Both the proposed FE-FV formulation in the current work as well as the traditional discrete formulation utilizing equal-order $P_1$ FE basis functions are implemented within the software framework BROOMStyx \cite{Sargado2019}, a general multiphysics simulator developed by the first author that allows for the combination of different discretization schemes with minimal overhead. Furthermore, the implementations for both formulations are optimized by storing quantities such as gradient matrices that do not vary with time at each evaluation point. Note that the incurred memory footprint of this is larger for the FE formulation in the phase-field subsystem, since it requires storage of both a 2-by-3 gradient matrix and 3-by-3 mass matrix. In contrast, the FV scheme only needs to store transmissibility coefficients which amount to a single scalar value per local cell edge. The code is designed to run in parallel on shared memory architectures by means of OpenMP directives, and dense matrix operations are carried out using BLAS and LAPACK routines provided by the Intel Math Kernel Library (MKL). Likewise, shared memory-parallel solution of sparse linear systems is obtained via the direct solver PARDISO that is included in MKL. We do not perform adaptive remeshing during the course of our simulations, hence the initial sparsity profiles of global coefficient matrix do not change. As in \cite{Sargado2018}, we exploit this property to speed up the solution of linear systems, specifically by performing symbolic factorization only at the beginning iteration of a substep, and thereafter proceeding directly to the numerical factorization phase of the solver in subsequent iterations. Additionally since both formulations result in symmetric global coefficient matrices, only the upper triangular portion of said matrices are actually assembled and manipulated by the linear solver which cuts down further the time for each iteration.
% ===========================================================================================
\section{Numerical examples}
We now compare performance of the proposed FE-FV formulation against an equal-order ($P_1$) FE approximation of the linear momentum and phase-field equations. For problems in 2D, plane strain behavior is assumed. Furthermore to put reported run times in context, we note that all simulations were carried out on a desktop machine equipped with a single Intel i7-8700K processor having 6 cores at 3.70 GHz base frequency. Each simulation run uses all 12 hyper-threads during execution of parallel regions such the assembly of global coefficient matrices and right hand sides, and also during solution of sparse linear systems. Finally, all meshes for 2D domains were generated using the software Gmsh \cite{Geuzaine2009}.

\subsection{Stationary crack in one dimension} \label{sec:stationary1D}
\begin{figure}
	\centering
	\includegraphics[width=0.9\textwidth]{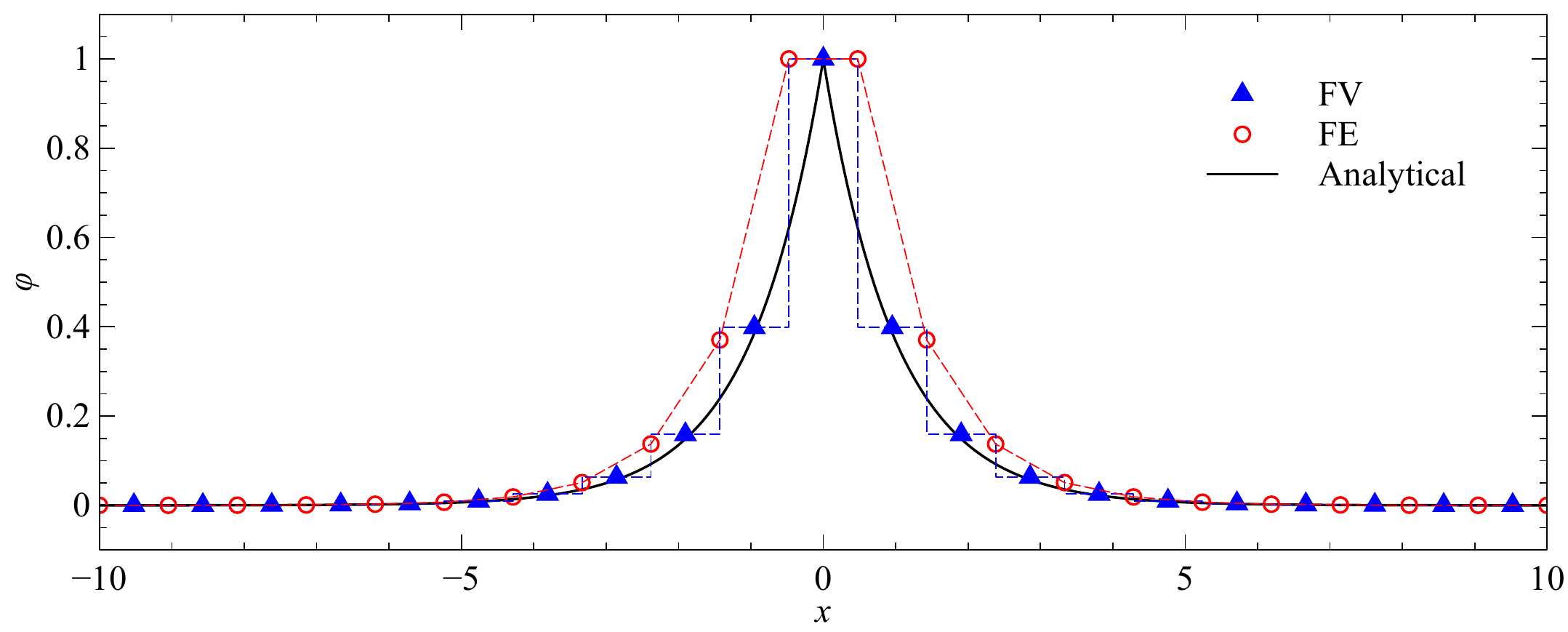}
	\caption{Numerical solutions for the phase-field corresponding to a fully developed crack at $x = 0$  obtained using cell-centered FVM (TPFA) versus FEM ($P_1$). Dashed lines show how $\phi$ is interpolated locally within each element.}
	\label{fig:comparison_1D}
\end{figure}

\begin{figure}
	\centering
	\includegraphics[width=0.45\textwidth]{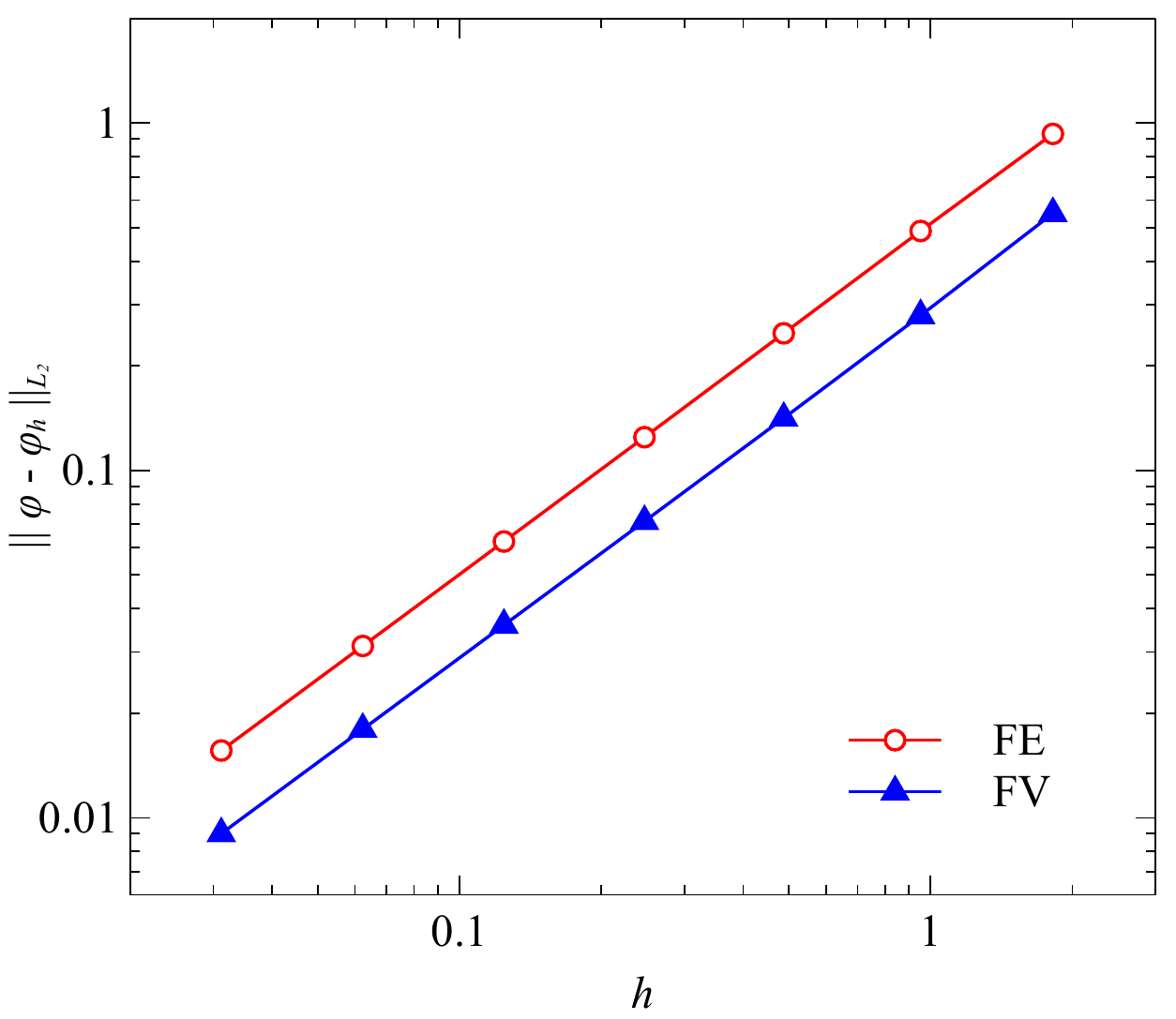}
	\caption{Convergence of FE and FV solutions of the phase-field profile with respect to mesh refinement for the 1-dimensional problem.}
	\label{fig:error_norms}
\end{figure}

In this example, we investigate the accuracy of using cell-centered finite volumes to solve the boundary value problem corresponding to a uniform cylindrical bar with endpoints at $x = \pm 10$ and which is completely cut by a crack at $x = 0$. The bar is unstressed, and we are concerned only with solving for the phase-field profile. Hence the originally coupled problem simplifies to the 1-dimensional ODE
\begin{linenomath}
\begin{equation}
	\ell \phi^{\prime\prime} \left( x \right) - \frac{1}{\ell} \phi \left( x \right) = 0 \quad \forall x \in \left( -10, 10 \right)
\end{equation}
\end{linenomath}
that is subject to the boundary condition $\phi^\prime \left( \pm 10 \right) = 0$, where for simplicity we set $\ell = 1$. The analytical solution in this case is given by
\begin{linenomath}
\begin{equation}
	\phi \left( x \right) = \exp \left( -\left| x \right| \right).
\end{equation}
\end{linenomath}
In a general discrete setting where the mechanics and phase-field equations are fully coupled, evolution of $\phi$ is driven by the source term $g^\prime \left( \phi \right) \psi^+ \left( \bm{\upepsilon} \right)$ so that the phase-field should achieve its maximum values at the discrete locations where the source term is applied. In a finite element approximation of the mechanics, stresses and strains are naturally defined at element Gauss points, which are then the natural choice of locations for application of the phase-field source term. In the case of $P_1$ FE, said Gauss points coincide with element midpoints. In order to achieve the same effect in the decoupled phase-field equation, the domain is discretized such that the coordinate $x = 0$ is located at the center of an element. In the FV simulation, the effect of a source term acting on this central element is mimicked by setting the value of its associated DOF to unity. On the other hand for the FE setup, the constraint $\phi = 1$ is applied at the nodes of the same element. Figure \ref{fig:comparison_1D} shows the results of both methods for a domain partitioned into 21 uniformly sized cells of size $h = 0.952$.
We can observe that although the phase-field is assumed to be constant within each cell for the FV solution, the value of $\phi$ at the center of each cell is much closer to the analytical solution compared to that obtained using finite elements. This can be attributed to the fact that approximation of the gradient terms at element endpoints by means of two-point fluxes effectively treats $\phi$ as discontinuous inside each cell, thus accommodating the presence of a cusp at $x = 0$. We also study the effect of mesh refinement on the accuracy of both methods by means of the $L_2$-norm of the error, defined as
\begin{linenomath}
\begin{equation}
	\| \phi - \phi_h \|_{L_2} = \sqrt{ \int_\Omega \left( \phi - \phi_h \right)^2 \dee\Omega}
\end{equation}
\end{linenomath}
Results are plotted in Figure \ref{fig:error_norms}, where $h \in \left\{ \frac{20}{11}, \frac{20}{21}, \frac{20}{41}, \frac{20}{81}, \frac{20}{161}, \frac{20}{321}, \frac{20}{641} \right\}$.
We can see that while the error converges at the same rate in the $L_2$-norm for both methods, the accuracy obtained with finite volumes is better, being about the same as that for a finite element solution using a mesh with element sizes reduced by two. Cell-centered finite volumes are thus a cost-effective alternative to $P_1$ FE for the fracture phase-field sub-problem, and in particular the factor two is found to hold also for higher dimensions as demonstrated in the succeeding examples.

% -------------------------------------------------------------------------------------------------------
\subsection{Tension test of notched specimen} \label{sec:tensionTest}
\begin{figure}
	\centering
	\begin{subfigure}[b]{0.35\textwidth}
		\centering
		\includegraphics[width=\textwidth]{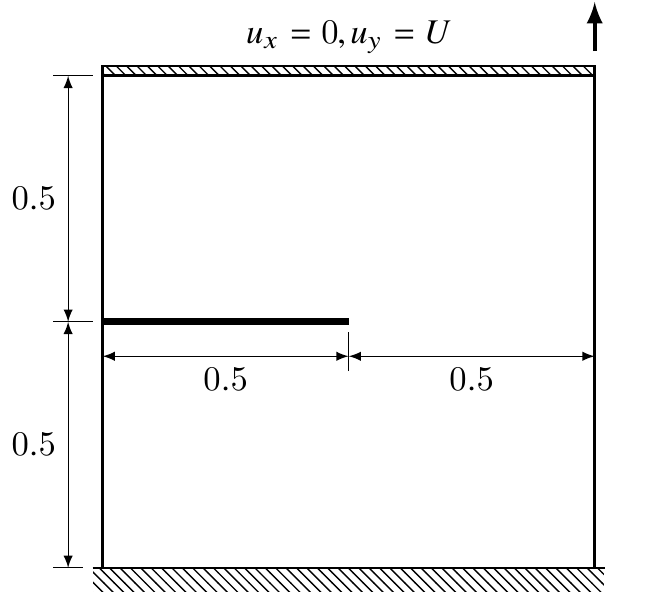}
		\caption{}
	\end{subfigure}
	\begin{subfigure}[b]{0.28\textwidth}
		\centering
		\includegraphics[width=\textwidth]{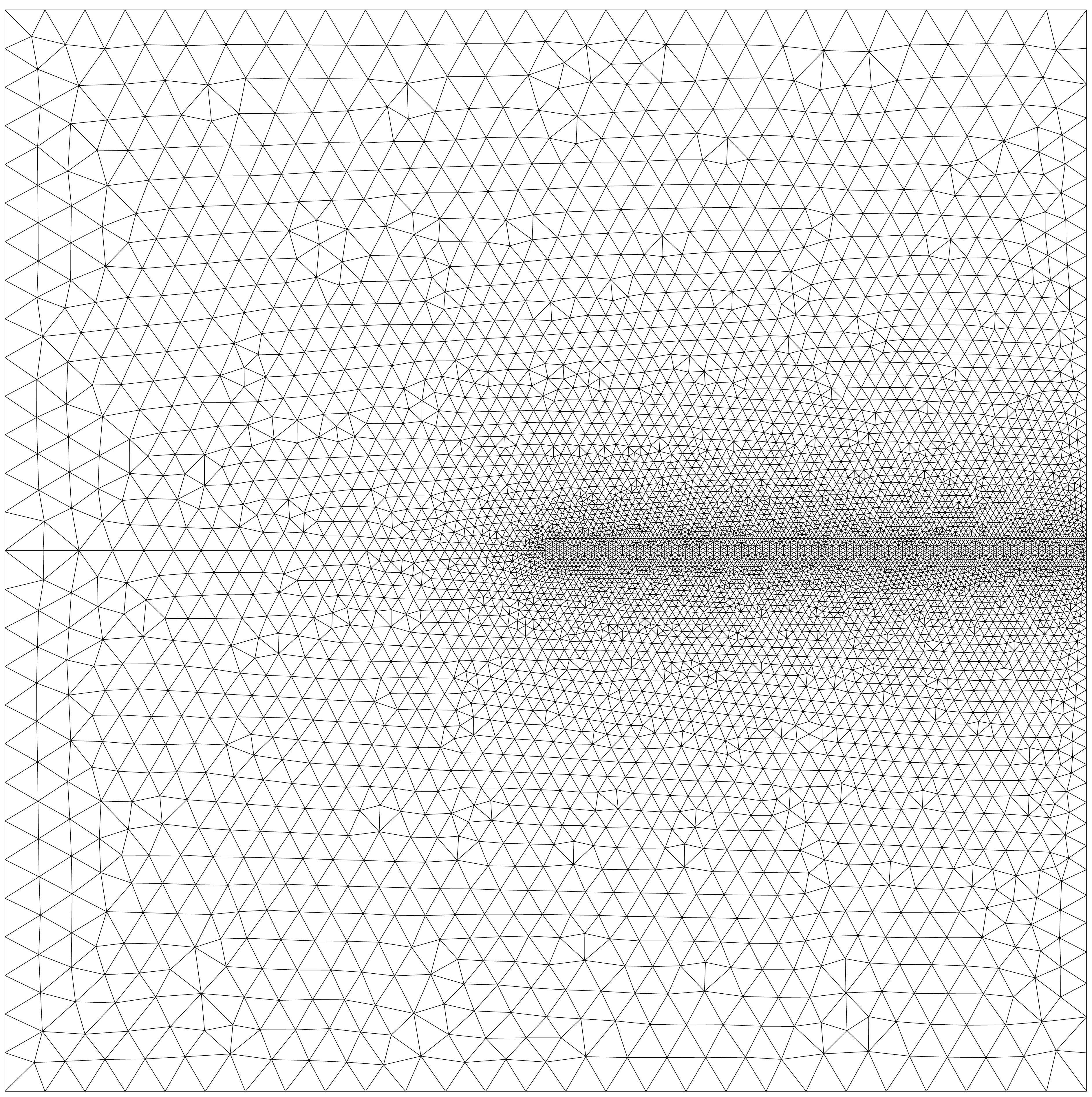}
		\caption{}
	\end{subfigure}
	\caption{Tension test of notched square specimen: (a) geometry and boundary conditions, (b) mesh refinement detail for $\ell/h =2$.}
	\label{fig:MieheTensionGeometry}
\end{figure}

The tensile (mode-I) fracturing of a notched square specimen is a benchmark problem originally introduced by Miehe et al. \cite{Miehe2010_ijnme} for demonstrating the capability of the phase-field approach to model brutal cracking. The specimen geometry and boundary conditions along with details of the a priori mesh refinement are shown in Figure \ref{fig:MieheTensionGeometry}, where the dimensions are given in mm.
Material parameters are set to the following: $E = 210,000$ MPa, $\nu = 0.3$, $\mathcal{G}_c = 2.7$ N/mm and $\ell = 0.0075$. The loading consists of a uniform vertical displacement $U = 0.00575$ mm imposed at the upper boundary of the specimen and initially applied in increments of $\Delta u_y = 1.0 \times 10^{-4}$ mm until $u_y = 0.00525$, and then in smaller increments of $\Delta u_y = 1.0 \times 10^{-5}$ up to a final displacement of $u_y = U$. As the resulting stresses are predominantly tensile, we use the isotropic model given in \eqref{eq:bourdinModel} for the bulk response. Two sets of simulations are conducted, one pertaining to standard FE-FE interpolation and another for the proposed FE-FV scheme. In each set, the setup described above is solved repeatedly over several mesh realizations corresponding different ratios between the phase-field length scale $\ell$ and the characteristic size $h$ of element edges at the refined region where the crack is expected to propagate. Note that for a given value of $\ell/h$, the same generated mesh is used in conjunction with both formulations.

\begin{figure}
	\centering
	\begin{minipage}{0.2\textwidth}
		\begin{subfigure}{\textwidth}
			\centering
			\includegraphics[width=\textwidth]{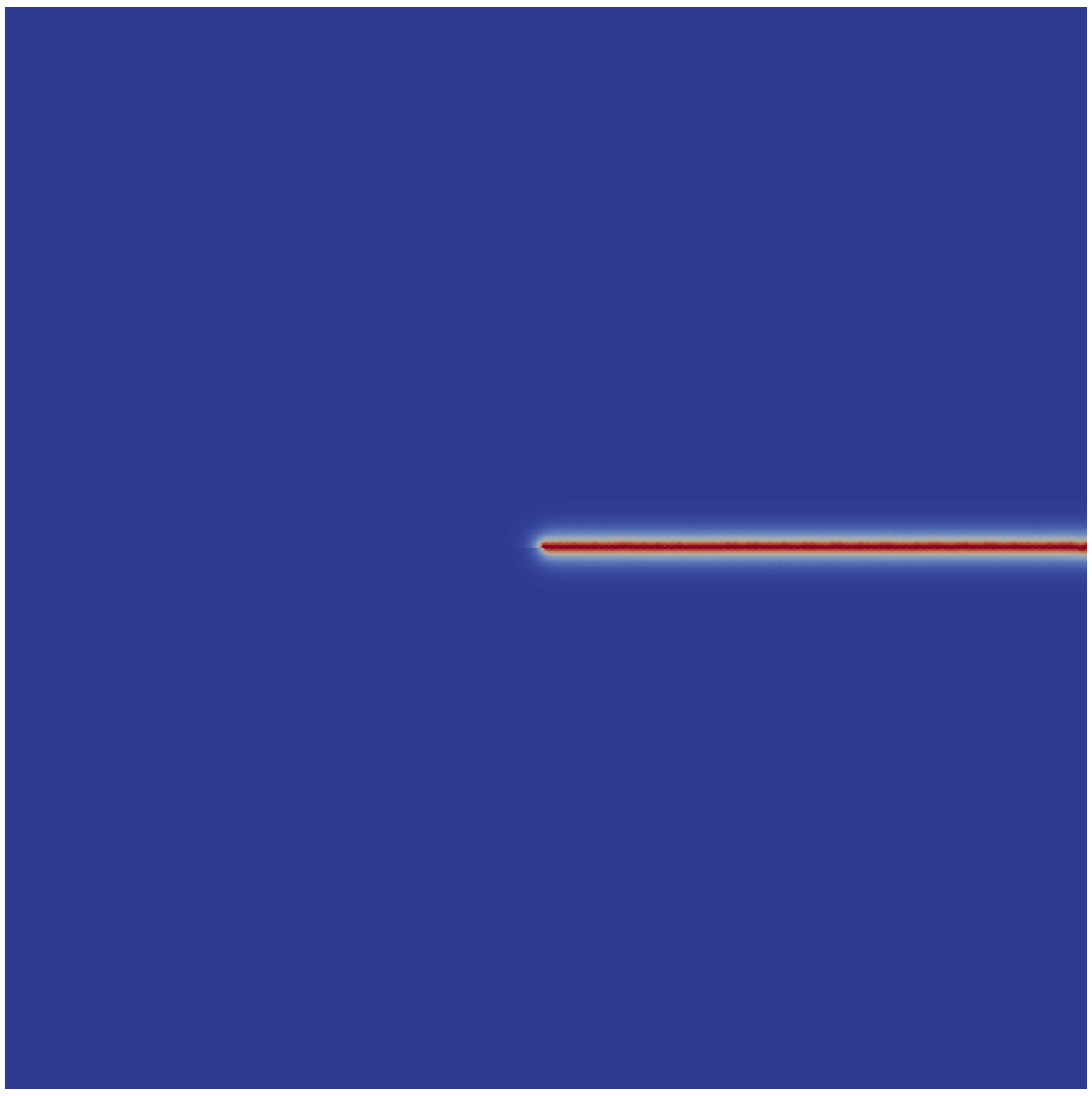}
			\caption{}
		\end{subfigure} \\
		\begin{subfigure}{\textwidth}
			\centering
			\includegraphics[width=\textwidth]{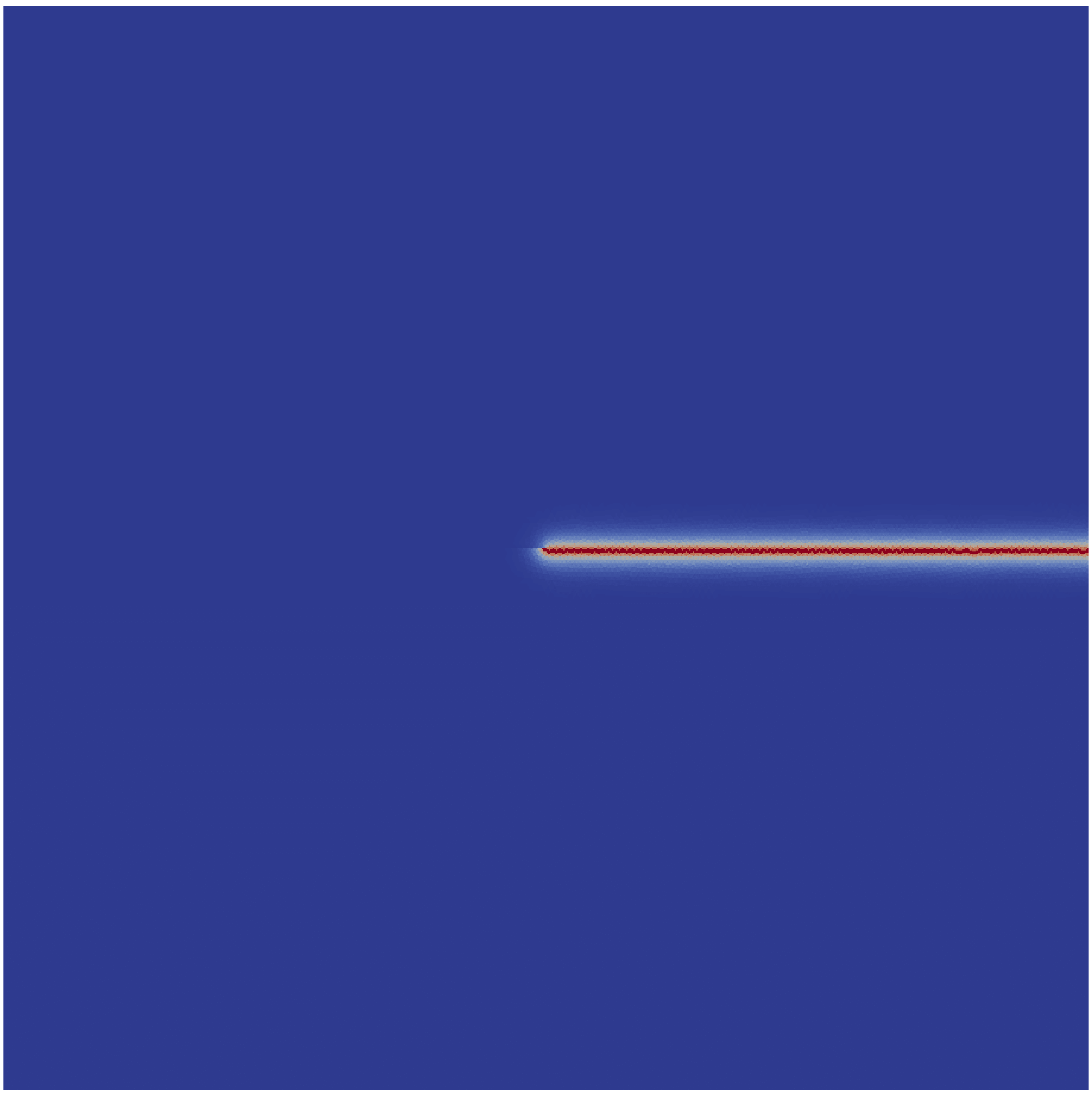}
			\caption{}
		\end{subfigure}
	\end{minipage} \hspace{5pt}
	\begin{minipage}{0.5\textwidth}
		\begin{subfigure}{\textwidth}
			\centering
			\includegraphics[width=\textwidth]{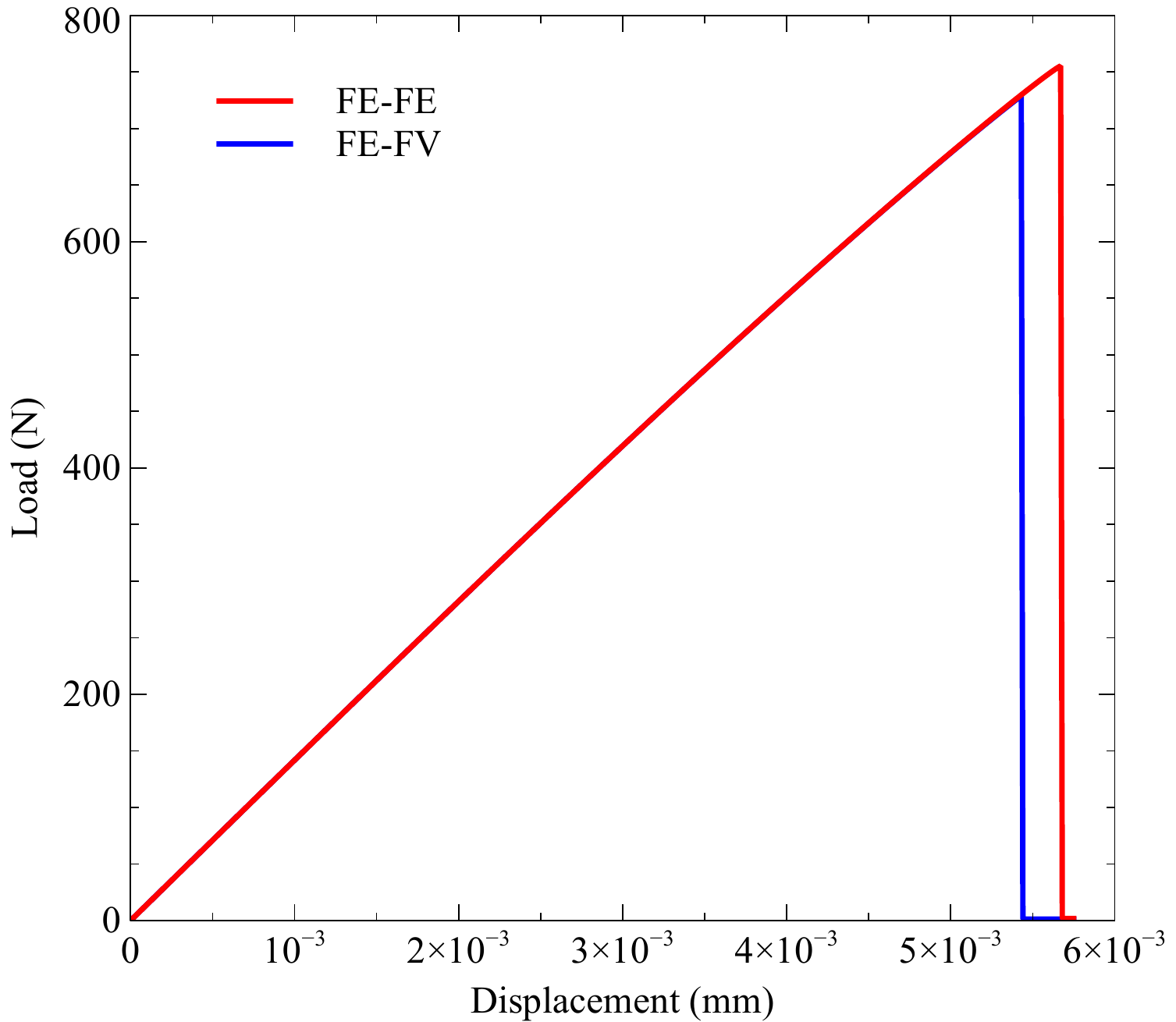}
			\caption{}
		\end{subfigure}
	\end{minipage}
	\caption{Results for notched tension test with $\ell/h = 2$. Final crack path obtained using (a) standard FE-FE, and (b) the FE-FV scheme. (c) Load displacement curves for the two schemes.}
	\label{eq:MieheTensionLoadDisp}
\end{figure}

\begin{figure}
	\centering
	\includegraphics[width=0.45\textwidth]{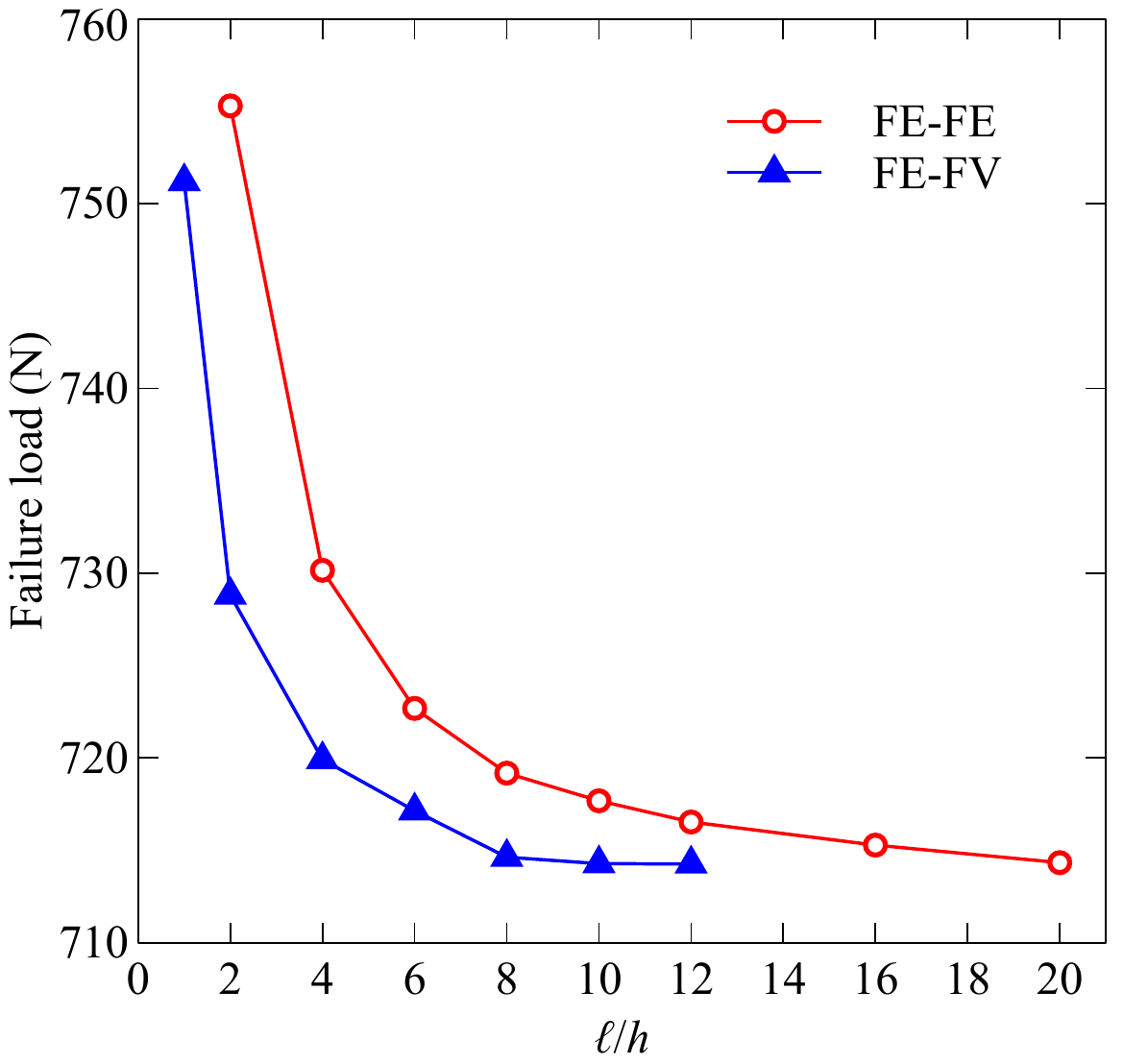}
	\caption{Failure loads for the notched tension test, obtained for different degrees of mesh refinement in conjunction with the proposed FE-FV scheme versus equal order FE interpolation.}
	\label{fig:failureLoadComparison}
\end{figure}

Table \ref{tab:tensionDetails} lists the resulting number of unknowns as well as the density of the global coefficient matrices for the simulation runs conducted.
\begin{table}
	\centering
	\caption{Details of individual simulations for the notched tension test.}
	\label{tab:tensionDetails}
	\begin{tabular}{crrrrrcccc}
		\toprule $\ell/h$ & \# Nodes & \# Cells & \# Unknowns, $\mathbf{u}$ & \multicolumn{2}{c}{\# Unknowns, $\phi$} & & \multicolumn{3}{c}{Global matrix density, \%} \\ \cmidrule{5-6} \cmidrule{8-10}
		& & & & \multicolumn{1}{c}{FE} & \multicolumn{1}{c}{FV} && $\mathbf{u}$ & $\phi$, FE & $\phi$, FV \\
		\midrule
		1 & 3,542 & 6,889 & 6,972\hspace{15pt} & ---\hspace{8pt} & 6,889 && 0.1060 & --- & 0.0361 \\
		2 & 6,567 & 12,906 & 13,022\hspace{15pt} & 6,567 & 12,906 && 0.0571 & 0.0604 & 0.0193 \\
		4 & 11,906 & 23,552 & 23,700\hspace{15pt} & 11,906 & 23,552 && 0.0315 & 0.0334 & 0.0106 \\
		6 & 17,031 & 33,784 & 33,950\hspace{15pt} & 17,031 & 33,784 && 0.0220 & 0.0234 & 0.0074 \\
		8 & 22,047 & 43,800 & 43,982\hspace{15pt} & 22,047 & 43,800 && 0.0170 & 0.0181 & 0.0057 \\
		10 & 27,832 & 55,360 & 55,552\hspace{15pt} & 27,832 & 55,360 && 0.0135 & 0.0143 & 0.0045 \\
		12 & 33,196 & 66,078 & 66,280\hspace{15pt} & 33,196 & 66,078 && 0.0113 & 0.0120 & 0.0038 \\
		16 & 43,297 & 86,266 & 86,482\hspace{15pt} & 43,297 & ---\hspace{8pt} &&  0.0087 & 0.0092 & --- \\
		20 & 53,776 & 107,212 & 107,440\hspace{15pt} & 53,776 & ---\hspace{8pt} && 0.0070 & 0.0074 & --- \\
		\bottomrule
	\end{tabular}
\end{table}
As both schemes use $P_1$ finite elements to discretize the linear momentum equation, for a given mesh they share the same number of displacement unknowns. The FE-FV formulation has around twice the number of unknowns for the phase-field subsystem as the traditional scheme, however the use of TPFA in former means that each row of the global coefficient matrix has at most four nonzero elements. This is substantially less connectivity than the corresponding finite element scheme as reflected in the last two columns of Table \ref{tab:tensionDetails}, where the global matrix density has been calculated as
\begin{linenomath}
\begin{equation}
	\text{Global matrix density, \%} = \frac{\text{\# nonzeros}}{\left( \text{\# unknowns} \right)^2} \times 100.
\end{equation}
\end{linenomath}
The final phase-field profiles corresponding to $\ell / h = 2$ are shown in Figure \ref{eq:MieheTensionLoadDisp} along with load-displacement curves for the two discretization schemes.
For the range of $\ell/h$ ratios tested in the current example, we found the resulting load displacement curves to be virtually identical for the two methods, the only observable difference being the actual instances of failure. The latter are summarized in Figure \ref{fig:failureLoadComparison}, where it can be seen that the critical load for the proposed FE-FV scheme has effectively converged with respect to mesh refinement at $\ell/h \approx 10$.
Moreover the results show that the failure load predicted using a mesh refinement of $\ell/h = n$ in conjunction with the proposed scheme is roughly equivalent to one obtained by utilizing mesh with $\ell/h = 2n$ for the FE-FE implementation, which corroborates the findings from Example \ref{sec:stationary1D}. We can thus infer that a critical load that is converged with respect to the mesh can be obtained with a mesh resolution of $\ell/h = 20$ in the crack path vicinity when using linear finite elements to discretize the phase-field. Incidentally, a mesh resolution of $\ell/h = 20$ was also used \cite{Giovanardi2017} to discretize the area around crack tips for a localized application of the phase-field model.

\begin{figure}
	\centering
	\begin{subfigure}{0.42\textwidth}
		\centering
		\includegraphics[width=\textwidth]{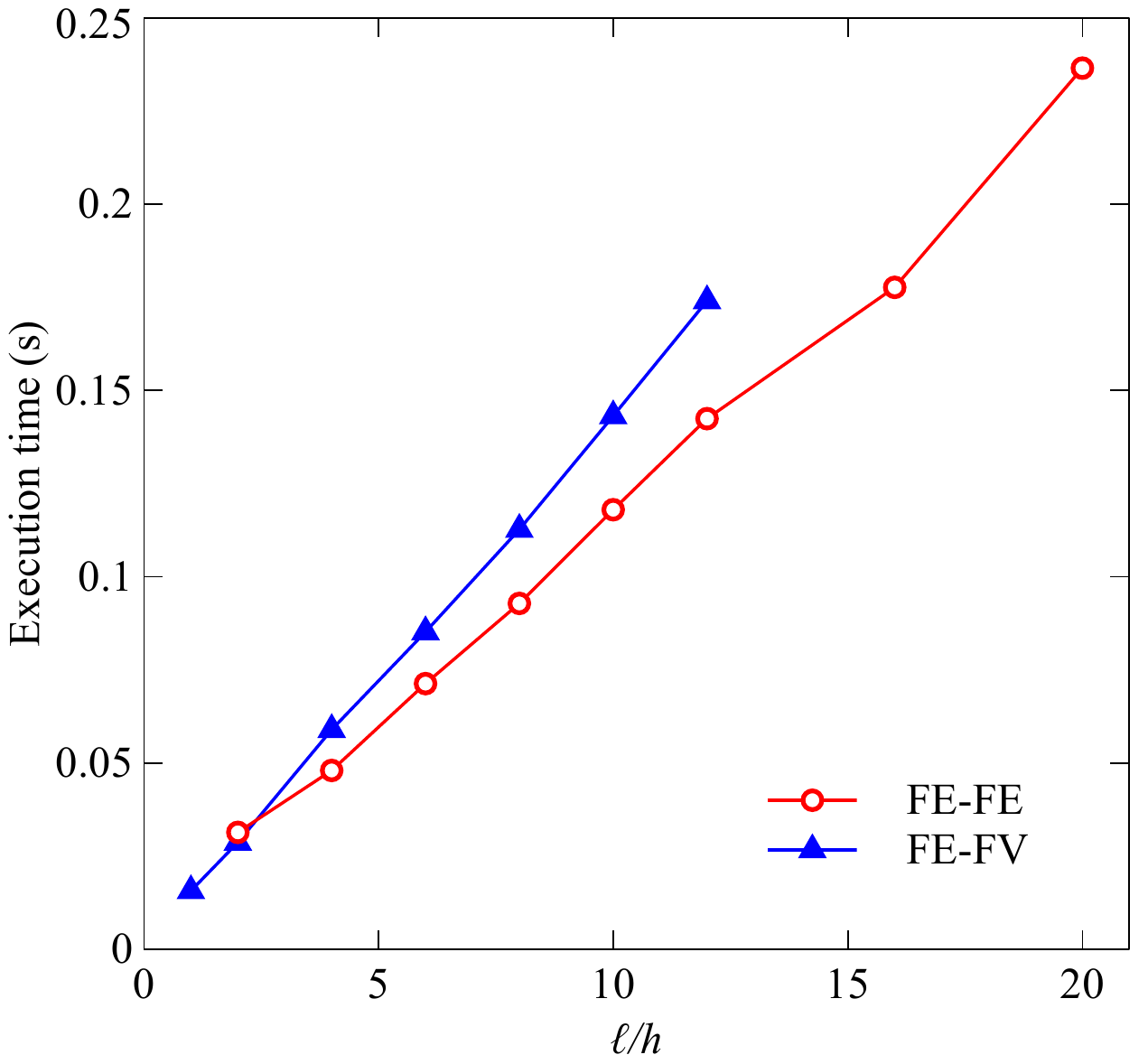}
		\caption{}
		\label{fig:iterTime}
	\end{subfigure} \hspace{10pt}
	\begin{subfigure}{0.4\textwidth}
		\centering
		\includegraphics[width=\textwidth]{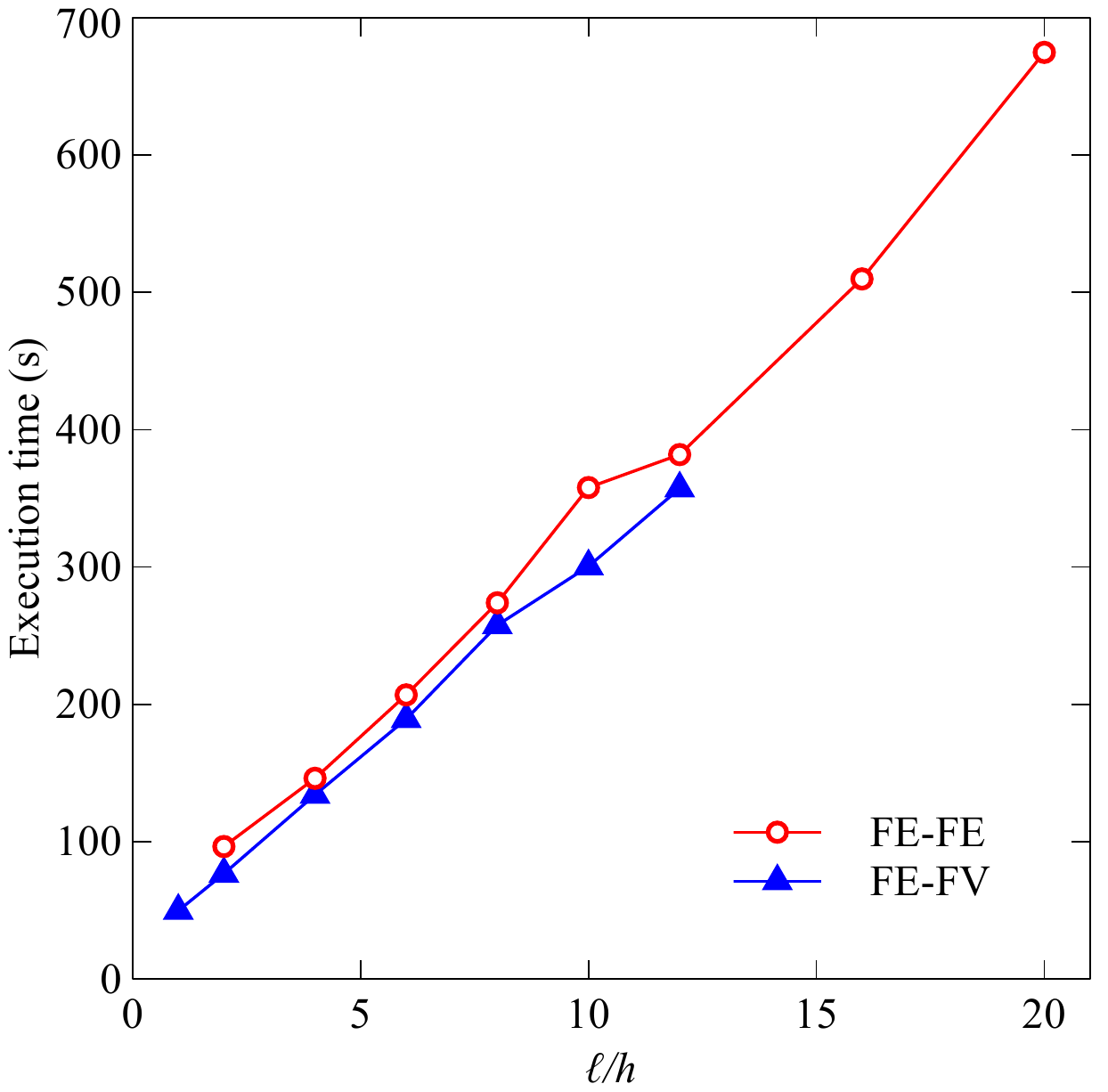}
		\caption{}
		\label{fig:totalTime}
	\end{subfigure}
	\caption{Comparison of execution times for the notched tension test: (a) typical run times for a single iteration of the alternate minimization algorithm; (b) total simulation times inclusive of writing output.}
\end{figure}

\begin{figure}
	\centering
	\includegraphics[width=0.42\textwidth]{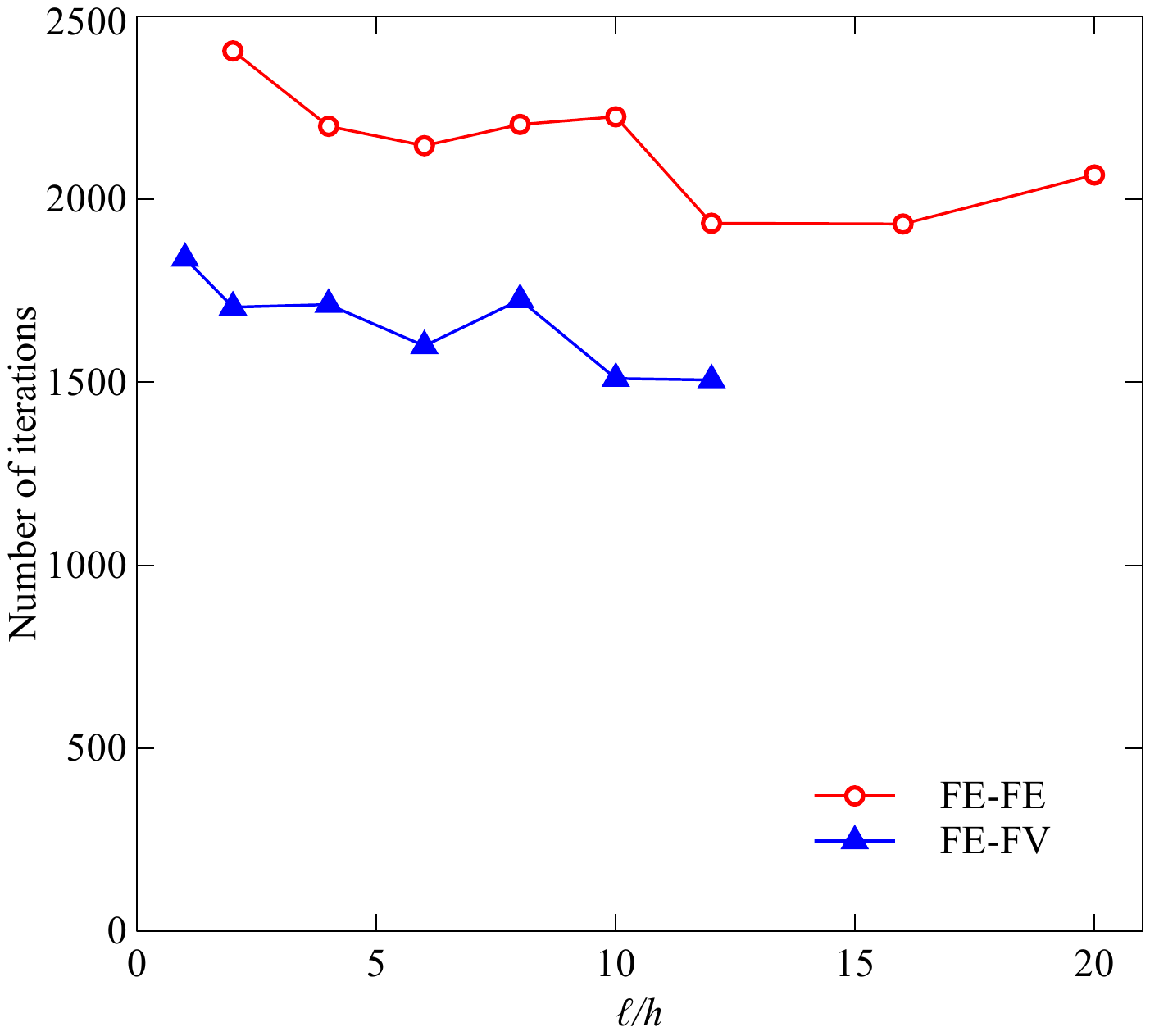}
	\caption{Variation of the total number of alternate minimization iterations with the choice of $\ell/h$ for the notched tension test.}
	\label{fig:totalIterations}
\end{figure}

On the other hand it is well known that setting $\ell/h < 2$ results in a mesh that is too coarse for the FE-FE scheme, leading to a severe overestimation of failure loads. Remarkably, not only does the FE-FV formulation allow for relatively coarse discretizations with $\ell/h = 1$, it even outperforms the standard formulation utilizing $\ell/h = 2$. This result becomes even more significant when we consider the overall run times for the two schemes. In Figure \ref{fig:iterTime}, we can observe that the time it takes to complete a single iteration of the alternate minimization algorithm is around 21\% more when using the combined discretization scheme as compared to the pure FE approach. 
This is somewhat to be expected, since the number of phase-field unknowns for the FE-FV scheme is roughly twice that of the FE-FE formulation. When we look at total simulation times however, the former turns out to be consistently cheaper (albeit only marginally) for all ratios of $\ell/h$ tested as shown in Figure \ref{fig:totalTime}, despite the fact that output files for the former are larger in size due to the increased number of unknowns. Such phenomenon can be explained upon examining the total number of iterations executed over the course of each simulation, which is plotted in Figure \ref{fig:totalIterations}.
Cumulative iterations are consistently less for the proposed formulation than for the traditional discretization by a difference of around 500 iterations, indicating that the alternate minimization algorithm converges faster for this problem with FE-FV scheme.

\begin{figure}
	\centering
	\includegraphics[width=0.45\textwidth]{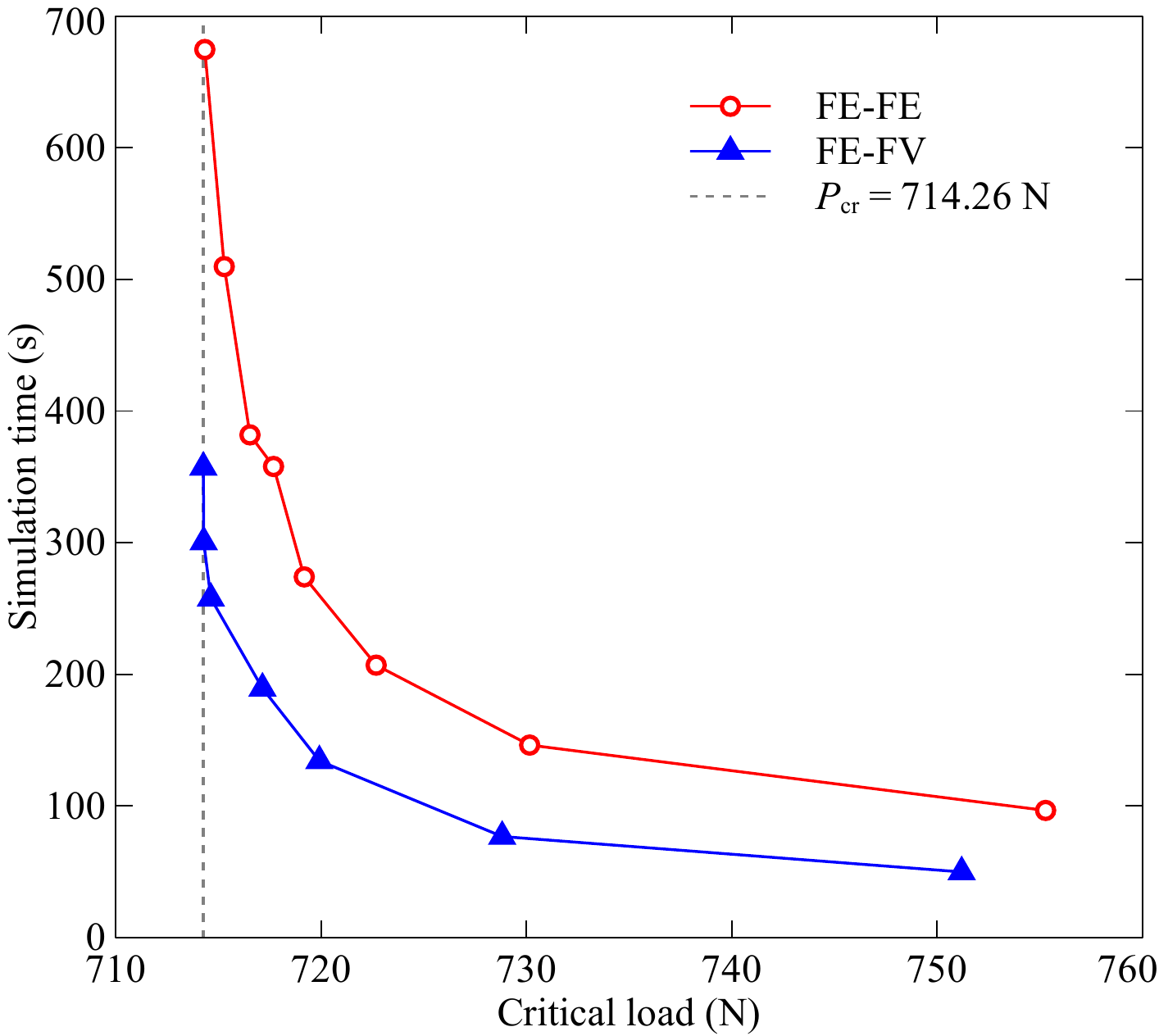}
	\caption{Plot of critical load vs.\ total simulation time for the notched tension test.}
	\label{fig:loadVsTime}
\end{figure}

An important point of comparison is the computational cost for obtaining a solution that is independent of the mesh for a given value of $\ell$. This can be determined by plotting failure loads against the total execution time, as shown in Figure \ref{fig:loadVsTime}.
Assuming that the FE-FE formulation also converges to the same critical load ($P_\text{cr} = 714.26$ N) with mesh refinement, we can infer that the computational cost of attaining mesh-insensitive results for a given $\ell$ is reduced by \emph{at least} 50\% when using the proposed FE-FV scheme. The aforementioned amount can be treated as representative of a lower bound due to the fact that we have confined the mesh refinement along crack path, which is known a priori for the given problem. If a uniform discretization is to be used (for instance when the crack trajectory is unknown at the beginning of the simulation), then the cost difference between the two formulations will be even greater than what is currently observed. Note that we make no conclusions regarding actual accuracy of the critical loads with respect to the true solution, since this requires that the proper value of $\ell$ be used, or that the model must be modified to account for usage of a smaller or larger $\ell$ \cite{Tanne2018,Sargado2018}.

% ------------------------------------------------------------------------------------------------------------------
\subsection{Shear loading of notched specimen}
\begin{figure}
	\centering
	\begin{subfigure}[b]{0.35\textwidth}
		\centering
		\includegraphics[width=\textwidth]{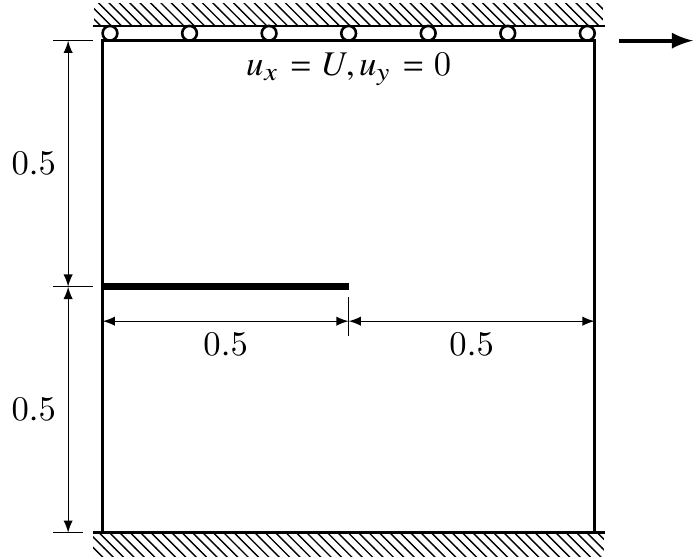}
		\caption{}
	\end{subfigure} \hspace{10pt}
	\begin{subfigure}[b]{0.28\textwidth}
		\centering
		\includegraphics[width=\textwidth]{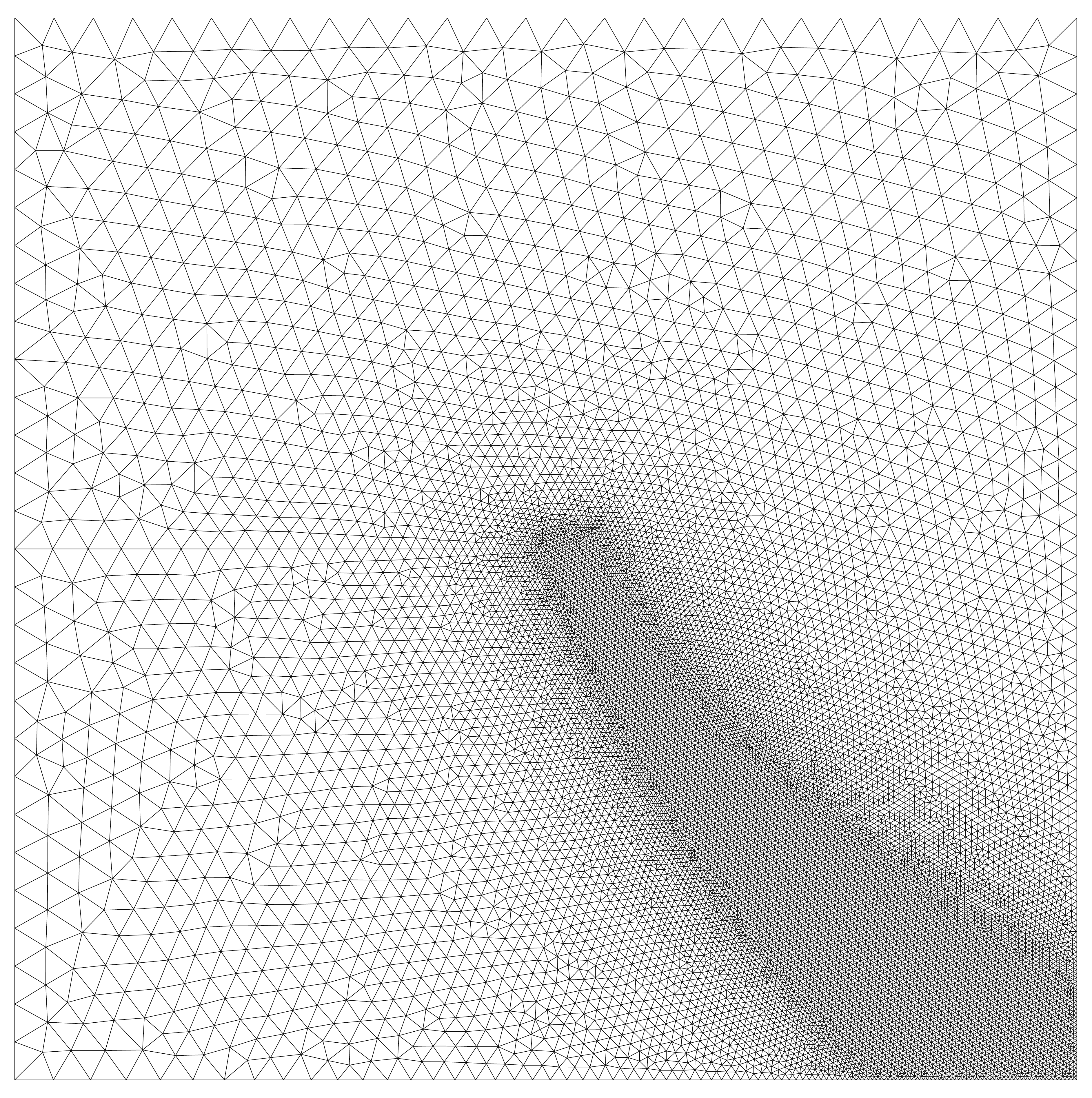}
		\caption{}
	\end{subfigure}
	\caption{Shear test of notched square specimen: (a) geometry and boundary conditions, and (b) mesh refinement detail for $\ell/h =2$.}
	\label{fig:MieheShearGeometry}
\end{figure}

The same notched specimen from the previous problem is now subjected to shear loading by prescribing a constant horizontal displacement at its top boundary as depicted in Figure \ref{fig:MieheShearGeometry}. As in the previous example, the mesh is refined a priori in the region where the crack is expected to propagate.
The present benchmark problem derives from \cite{Miehe2010_ijnme,Miehe2010_cmame}, however we note that there exists some variation across different studies on exactly what boundary conditions to apply. Our treatment here is similar to \cite{Ambati2015} in which the left and right sides of the specimen are kept traction free. Material parameters as well as the phase-field regularization parameter are the same as in the previous example, however for the current problem we make use of the anisotropic bulk degradation model of Amor et al.\ \cite{Amor2009} in order to prevent the formation of cracks under compression. The imposed horizontal displacement at the top boundary is applied in increments of $\Delta u_x = 1 \times 10^{-4}$ mm until a cumulative value of $u_x = 0.0092$, and thereafter the increments are decreased to $\Delta u_x = 1 \times 10^{-5}$ until the final displacement of $U = 0.015$ mm is reached. A total of five simulation runs are conducted: two each for $\ell/h = 2$ and $\ell = 4$, and an additional run corresponding to $\ell/h=1$ for the mixed FE-FV scheme. Details for individual simulations are listed in Table \ref{tab:mieheShearDetails}, wherein the solution times reported are from the beginning of each simulation until convergence of the time step corresponding to an applied displacement of $u_x = 0.01344$ mm at the top boundary of the specimen.
\begin{table}
	\centering
	\caption{Details of individual simulations for the notched shear test.}
	\label{tab:mieheShearDetails}
	\begin{tabular}{ccccccc}
		\toprule Run \# & Method & $\ell / h$ & \multicolumn{3}{c}{\# Unknowns} & Run time (hr : min : sec) \\ \cmidrule{4-6}
		&&& $\mathbf{u}$ & $\phi$ & Total & (until $u_x = 0.01344$ mm) \\
		\midrule 1 & FE-FE & 2 & 28,314 & 14,295 & {\color{white}{0}}42,609 & 01 : 23 : 44 \\
		2 & FE-FE & 4 & 83,676 & 42,042 & 125,718 & 06 : 04 : 48 \\
		3 & FE-FV & 1 & 11,078 & 11,028 & {\color{white}{0}}22,106 & 00 : 20 : 18 \\
		4 & FE-FV & 2 & 28,314 & 28,261 & {\color{white}{0}}56,575 & 01 : 11 : 43 \\
		5 & FE-FV & 4 & 83,676 & 83,632 & 167,308 & 06 : 42 : 59 \\
		\bottomrule
	\end{tabular}
\end{table}
The reason for choosing this particular endpoint for measuring run times is that some of the simulations fail to achieve convergence in the succeeding time steps; the reported run times are intended to be more representative of the ideal situation in which all substeps converge. We can see that for the finest discretization ($\ell/h = 4$), use of the FE-FV formulation results in substantially longer execution time than the standard FE-FE discretization. This is due to a higher number of total iterations performed in the former, which is contrary to behavior exhibited in the previous example wherein cumulative iteration counts for the FE-FV scheme were consistently fewer than for the traditional formulation. It is not yet clear at this point why we obtain such behavior, since for instance when $\ell/h = 2$ in the current example, the simulation utilizing the FE-FV discretization is faster. We note that for nonlinear problems, the number of iterations may depend heavily on the relative size of time steps and also the tolerance criteria chosen by the user, however an in-depth investigation of this is beyond the scope of the present work. We further note that while we have utilized an unmodified version of the alternate minimization algorithm as described in \cite{Bourdin2008} for the present study, a substantial amount of effort has been exerted towards reducing the number of iterations in the nonlinear solution algorithms used to solve the coupled mechanics/phase-field problem, e.g.\ \cite{Gerasimov2016,Farrell2017,Wick2017}.

\begin{figure}
	\centering
	\includegraphics[width=0.5\textwidth]{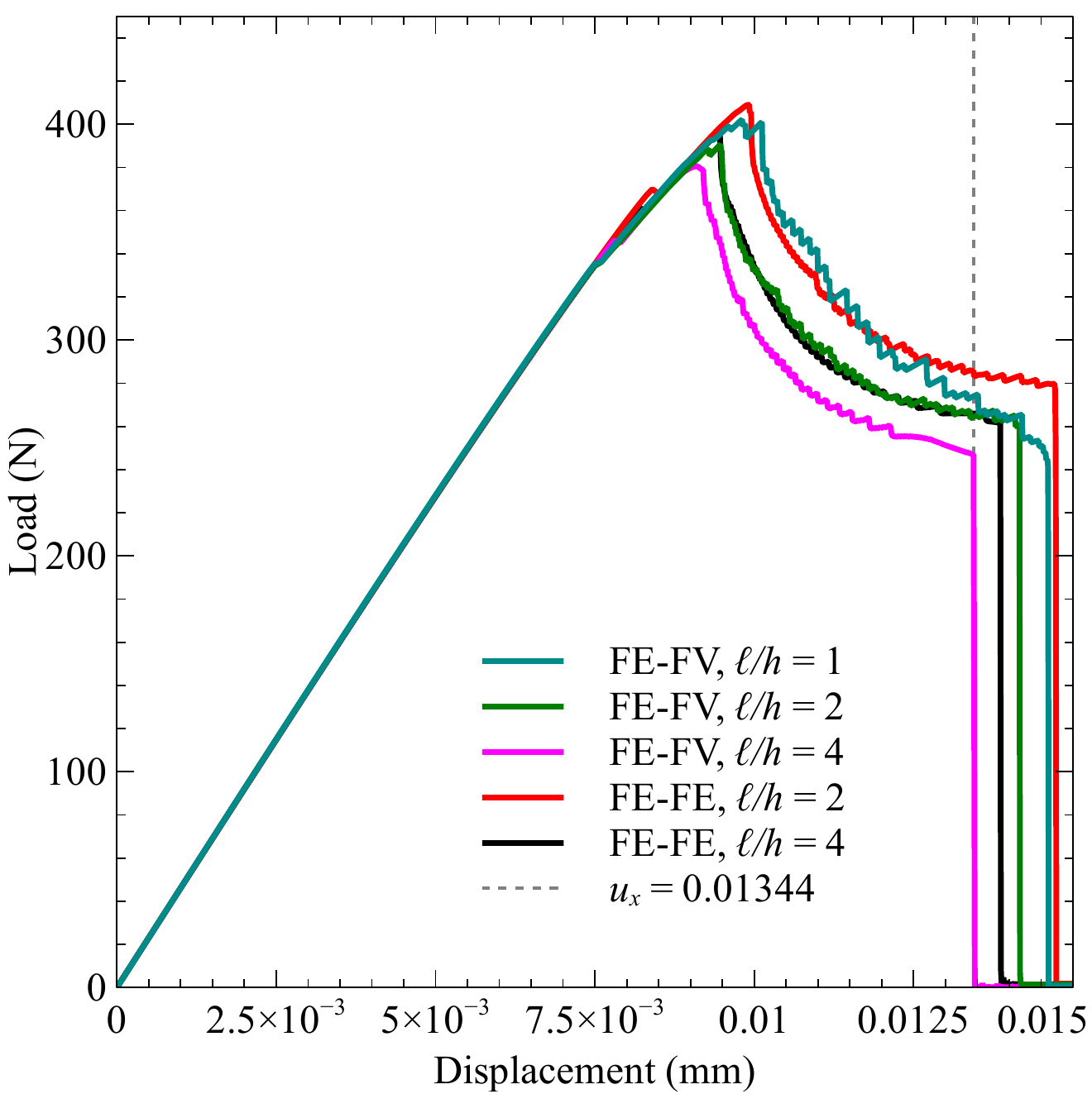}
	\caption{Superimposed load-displacement curves for the notched shear test.}
	\label{fig:shearLoadDisp}
\end{figure}

\begin{figure}
	\centering
	\begin{subfigure}{0.29\textwidth}
		\includegraphics[width=\textwidth]{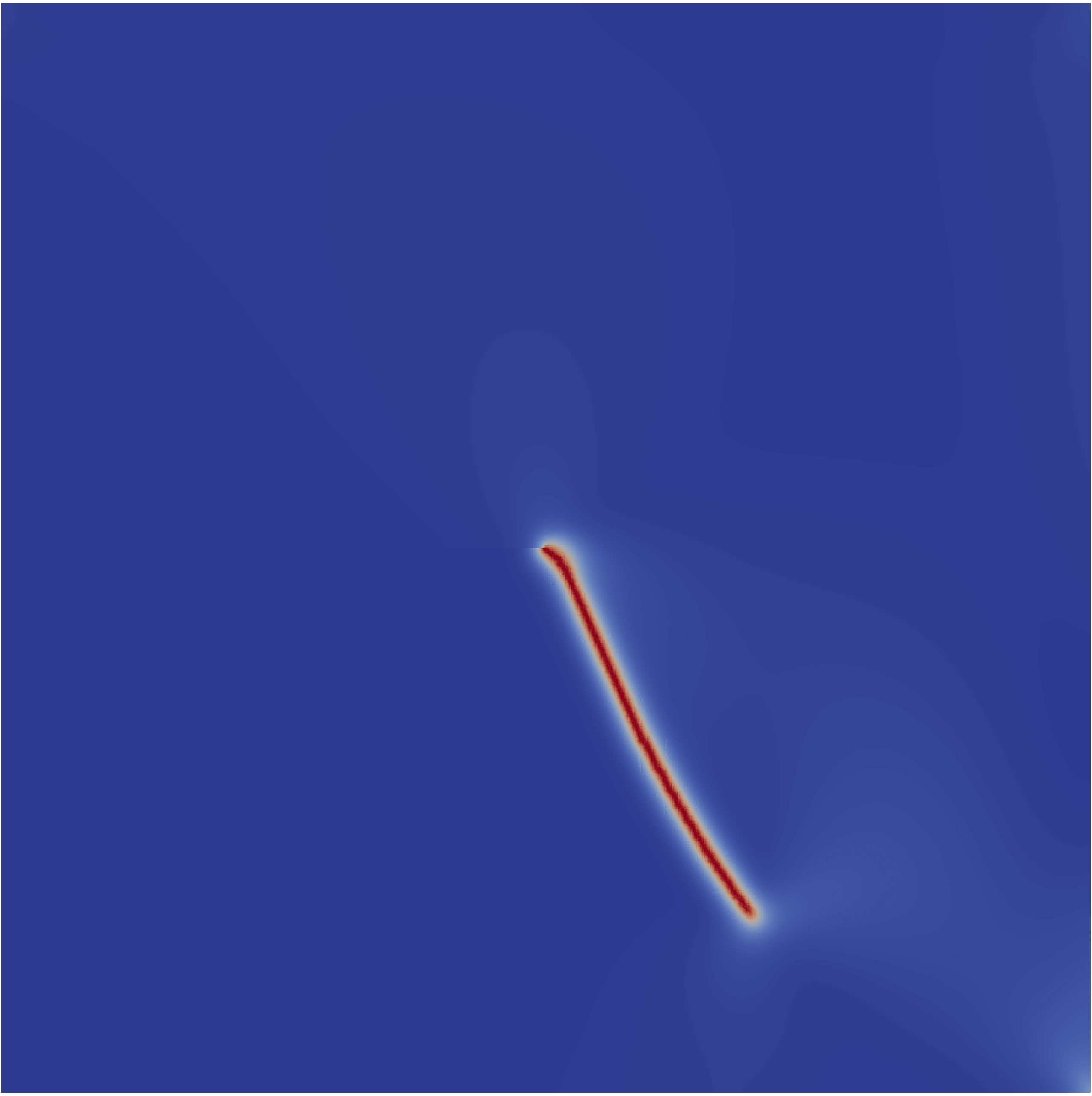}
		\caption{}
	\end{subfigure}
	\begin{subfigure}{0.29\textwidth}
		\includegraphics[width=\textwidth]{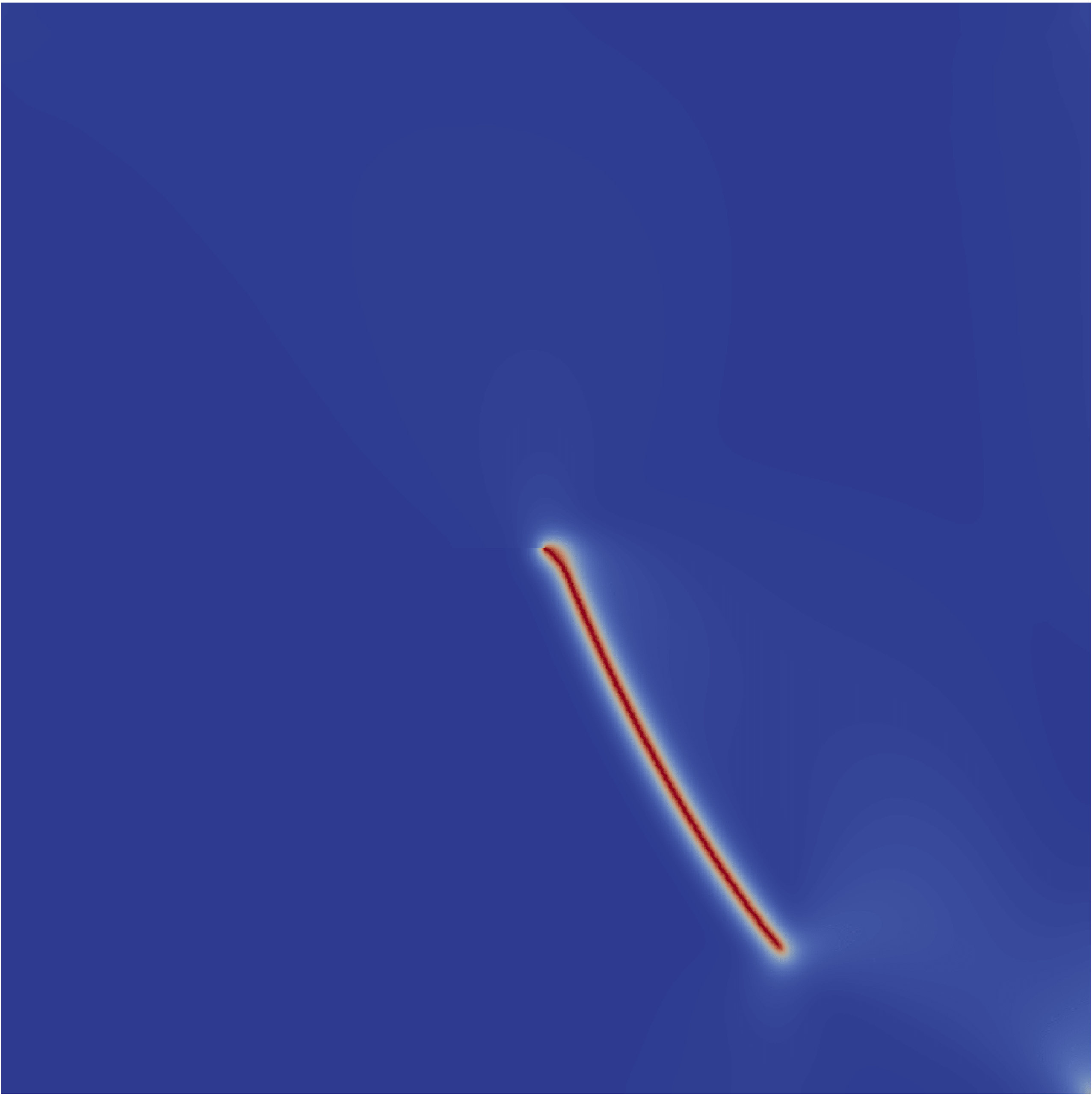}
		\caption{}
	\end{subfigure} \hspace*{0.29\textwidth} \\
	\begin{subfigure}{0.29\textwidth}
		\includegraphics[width=\textwidth]{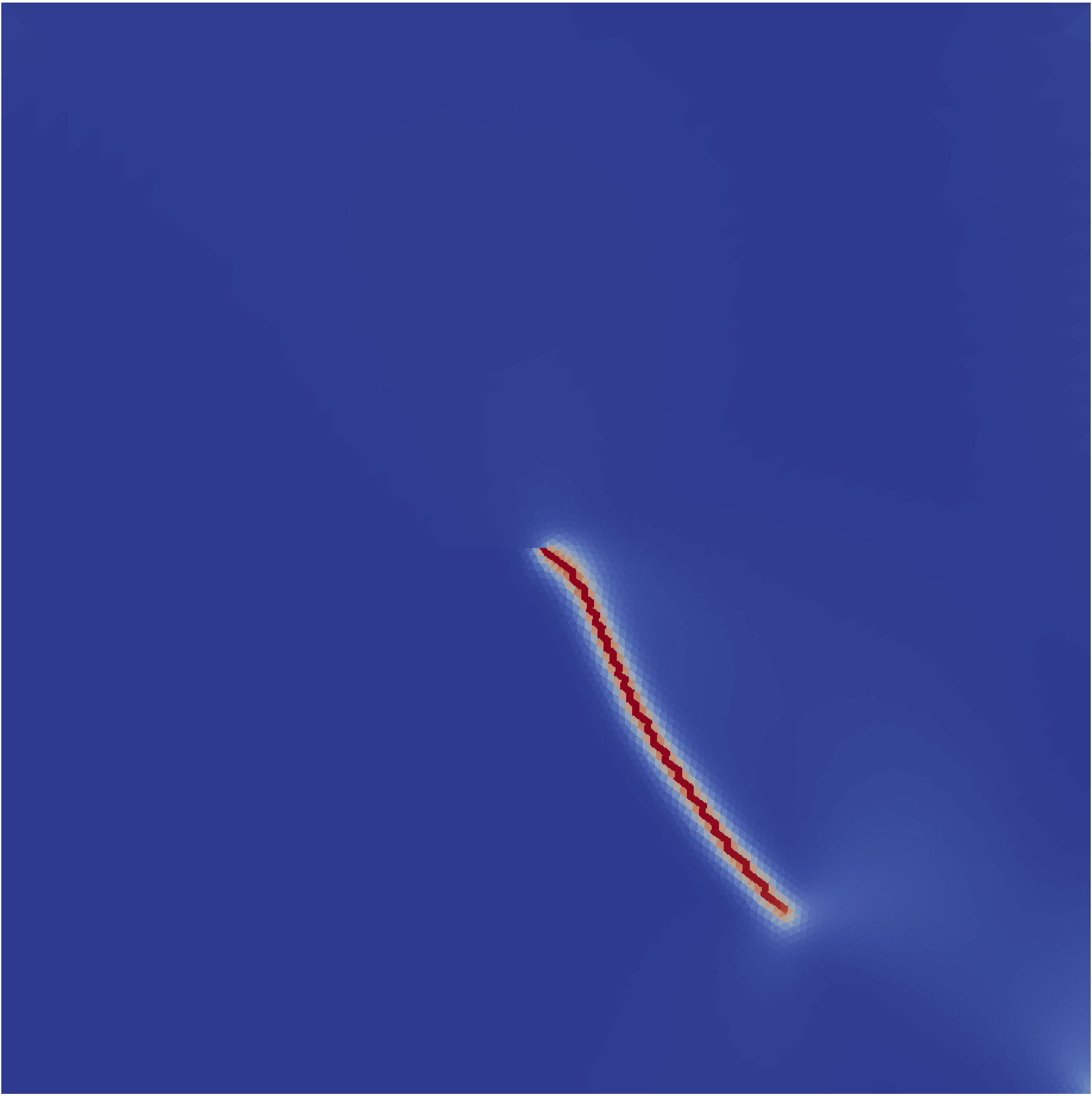}
		\caption{}
	\end{subfigure}
	\begin{subfigure}{0.29\textwidth}
		\includegraphics[width=\textwidth]{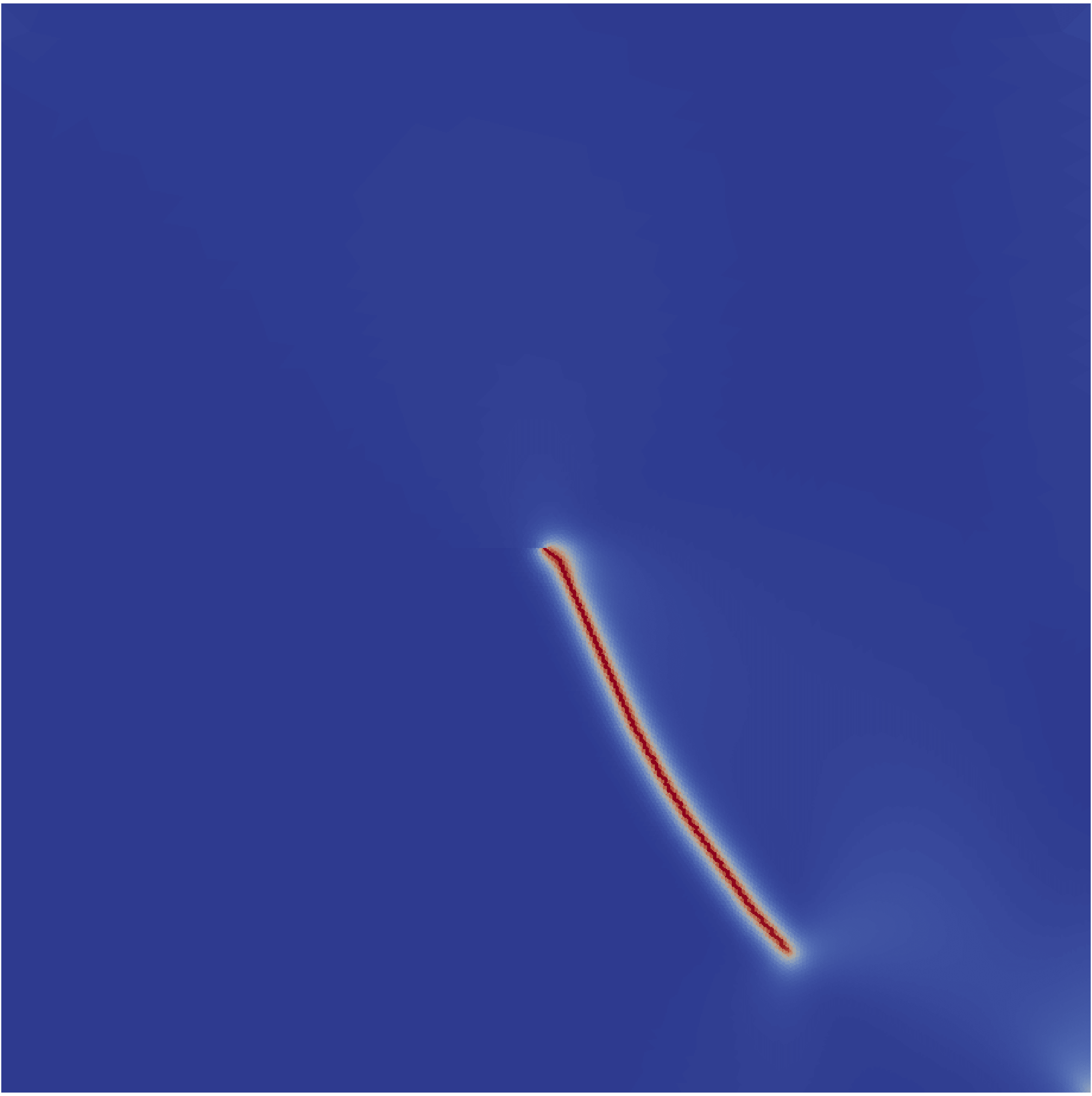}
		\caption{}
	\end{subfigure}
	\begin{subfigure}{0.29\textwidth}
		\includegraphics[width=\textwidth]{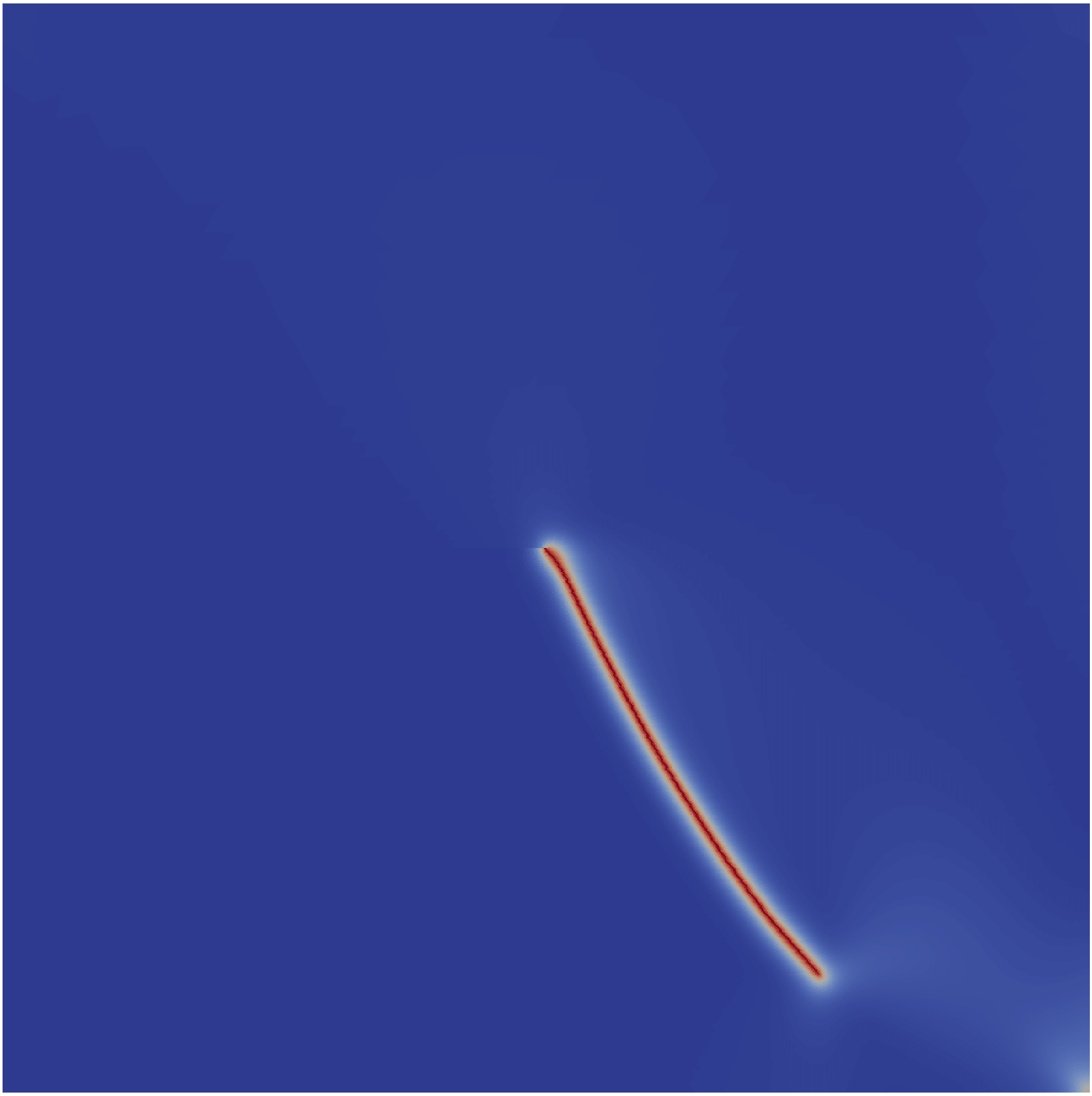}
		\caption{}
	\end{subfigure}
	\caption{Comparison of crack paths for the notched shear test corresponding to $u_x = 0.01344$: (a) FE-FE with $\ell/h=2$, (b) FE-FE with $\ell/h=4$, (c) FE-FV with $\ell/h=1$, (d) FE-FV with $\ell/h=2$, (e) FE-FV with $\ell/h=4$.}
	\label{fig:shearCrackPaths}
\end{figure}

Load-displacement curves obtained from the different simulations are plotted in Figure \ref{fig:shearLoadDisp} and exhibit the same qualitative behavior as reported in \cite{Ambati2015}, with the difference in predicted peak loads attributable to the fact that we use a slightly larger value of $\ell$ in our simulations.
More importantly, the results seem to confirm the trend shown in the prior example concerning the relationship between $\ell/h$ and the critical load, namely that FE-FE discretization of the governing equations requires a finer resolution of the mesh to yield a similar prediction of the mechanical response as the FE-FV formulation. In particular the post-peak behavior modeled by means of the former in conjunction with a mesh characterized by $\ell/h = 4$ and the latter with $\ell/h = 2$ are nearly identical. The same observation can be made by comparing the simulated crack paths shown in Figure \ref{fig:shearCrackPaths}.
Taking into account that the primary value of such simulations is in their predictive capability, the proper comparison that should be made with regard to computational cost is between runs 1 and 3, and similarly between runs 2 and 4 in Table \ref{tab:mieheShearDetails}. From the reported run times, we can see that adoption of the FE-FV scheme in place of the traditional FE-FE formulation results in a reduction of around 75\% in computational cost for the former case, and around 80\% for the latter.

\subsection{Stretching of specimen with multiple preexisting cracks} \label{sec:manyCracks}
\begin{figure}
	\centering
	\includegraphics[width=0.4\textwidth]{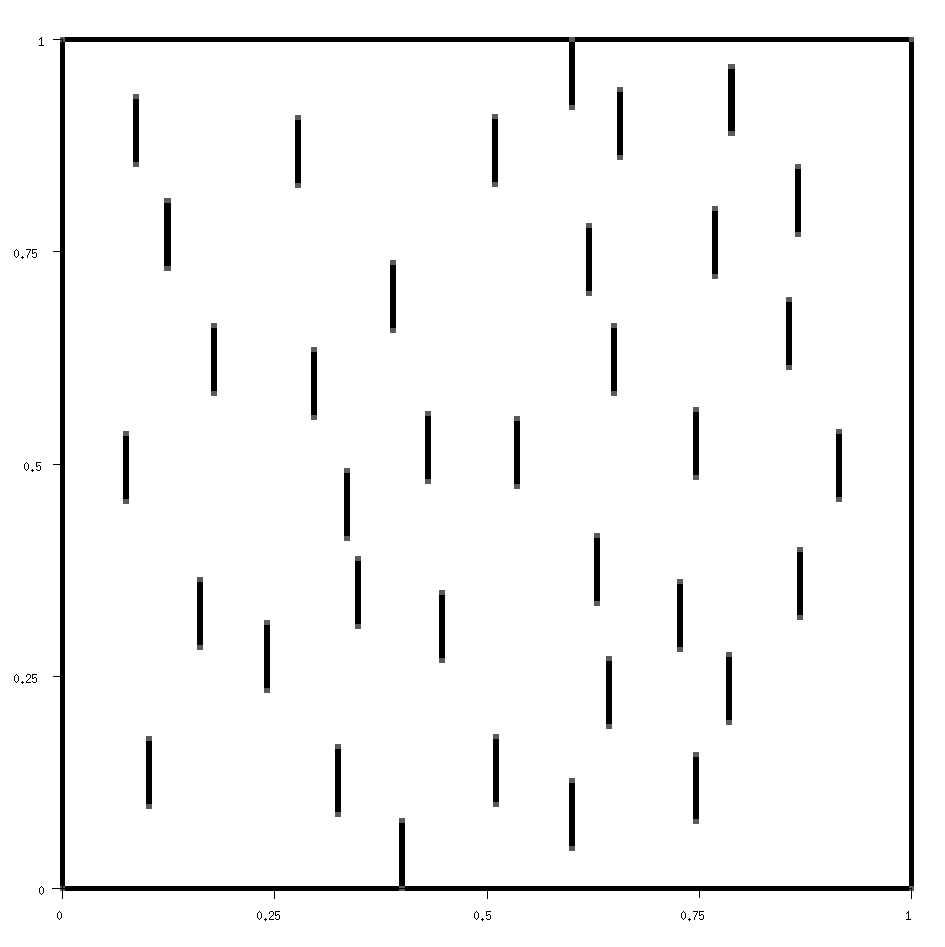}
	\caption{Geometry for specimen with preexisting cracks.}
	\label{fig:manyCracksGeometry}
\end{figure}

\begin{figure}
	\centering
	\includegraphics[width=0.5\textwidth]{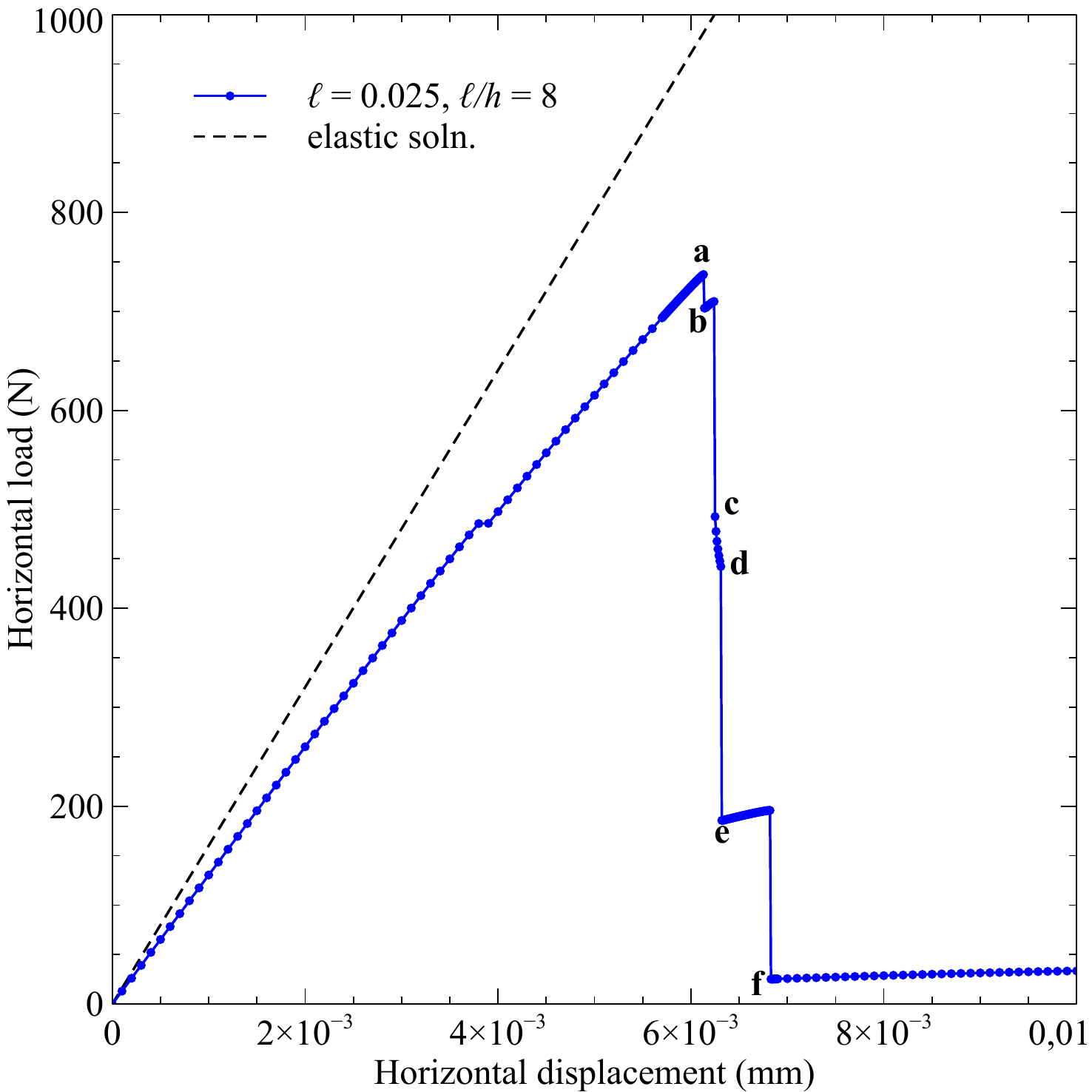}
	\caption{Load-displacement curve for the specimen of Example \ref{sec:manyCracks}, with elastic solution plotted for comparison. Dots correspond to solutions at actual time steps, while letter annotations indicate snapshots of the crack evolution shown in Figure \ref{fig:manyCrackEvolution}.}
	\label{fig:loadDisp_manyCracks}
\end{figure}

\begin{figure}
	\centering
	\begin{subfigure}{0.3\textwidth}
		\includegraphics[width=\textwidth]{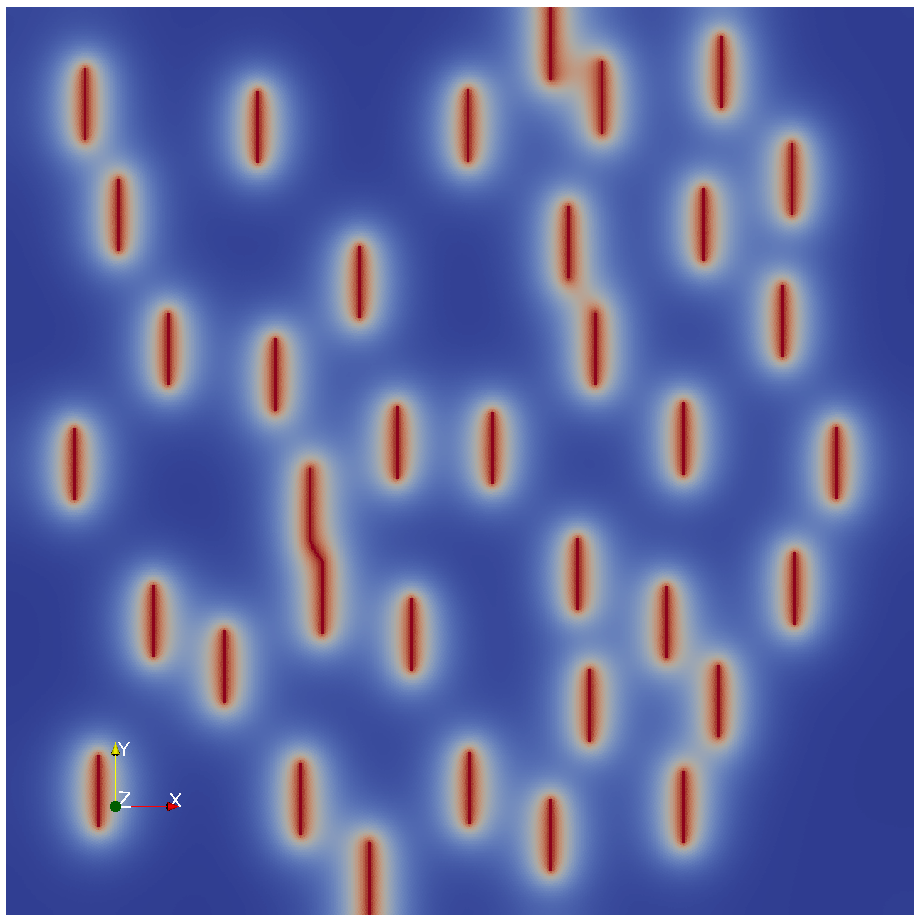}
		\caption{}
	\end{subfigure}
	\begin{subfigure}{0.3\textwidth}
		\includegraphics[width=\textwidth]{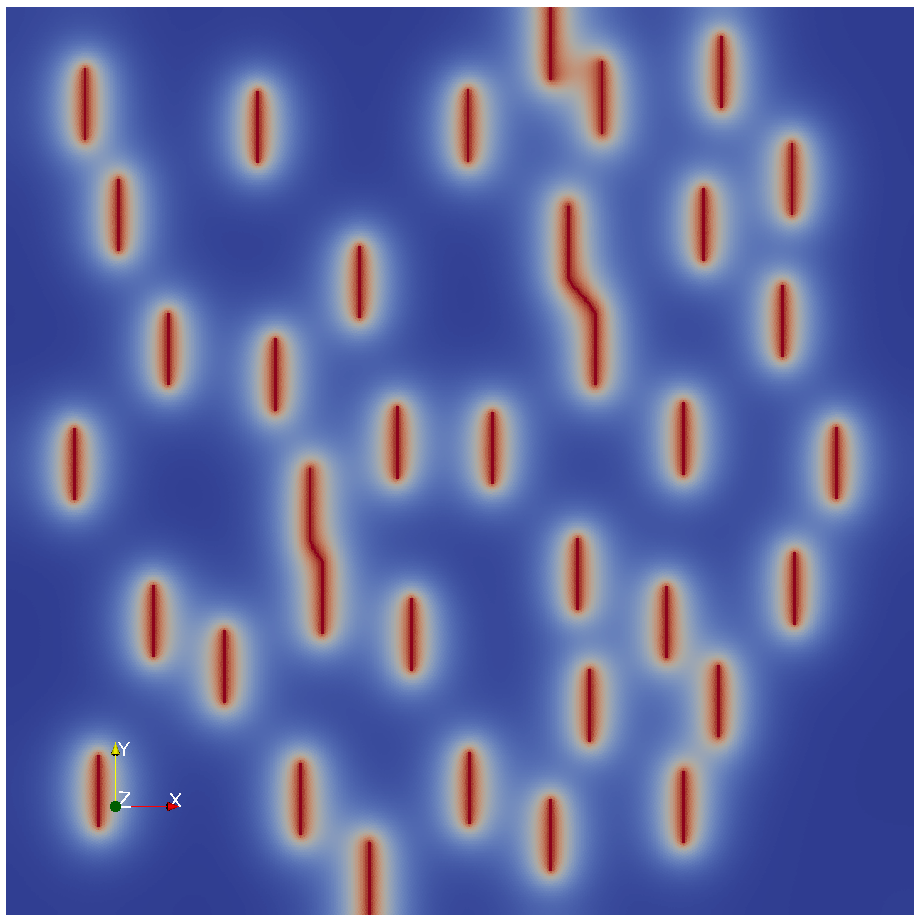}
		\caption{}
	\end{subfigure}
	\begin{subfigure}{0.3\textwidth}
		\includegraphics[width=\textwidth]{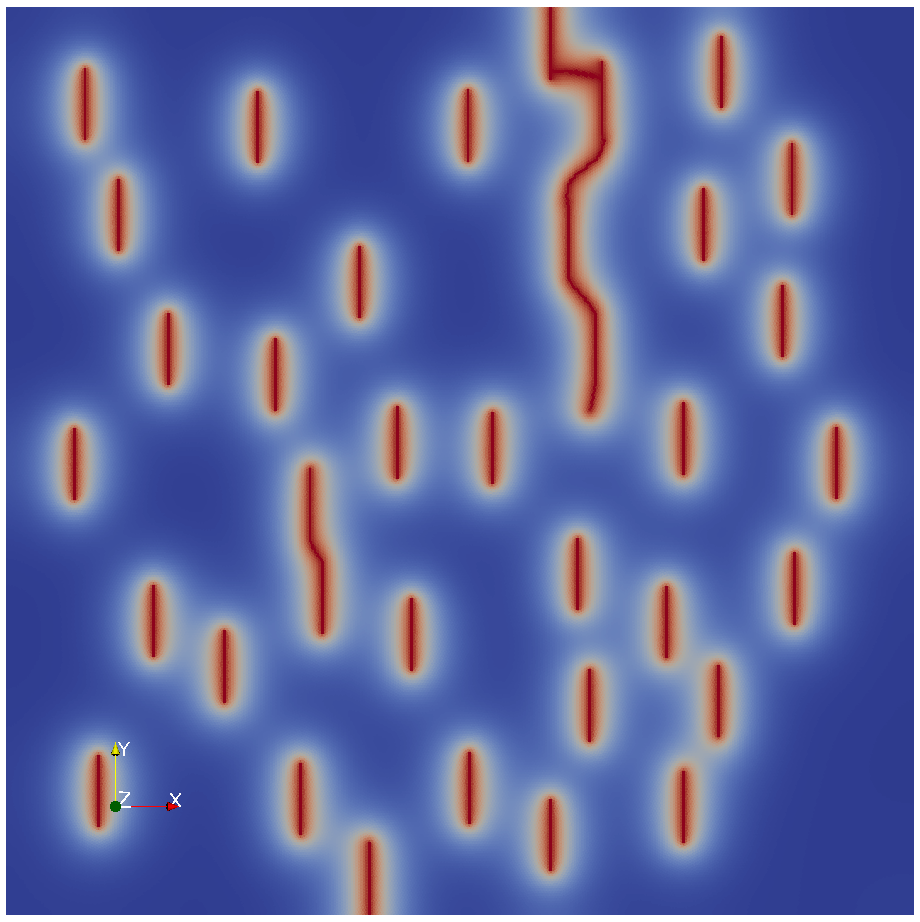}
		\caption{}
	\end{subfigure} \\
	\begin{subfigure}{0.3\textwidth}
		\includegraphics[width=\textwidth]{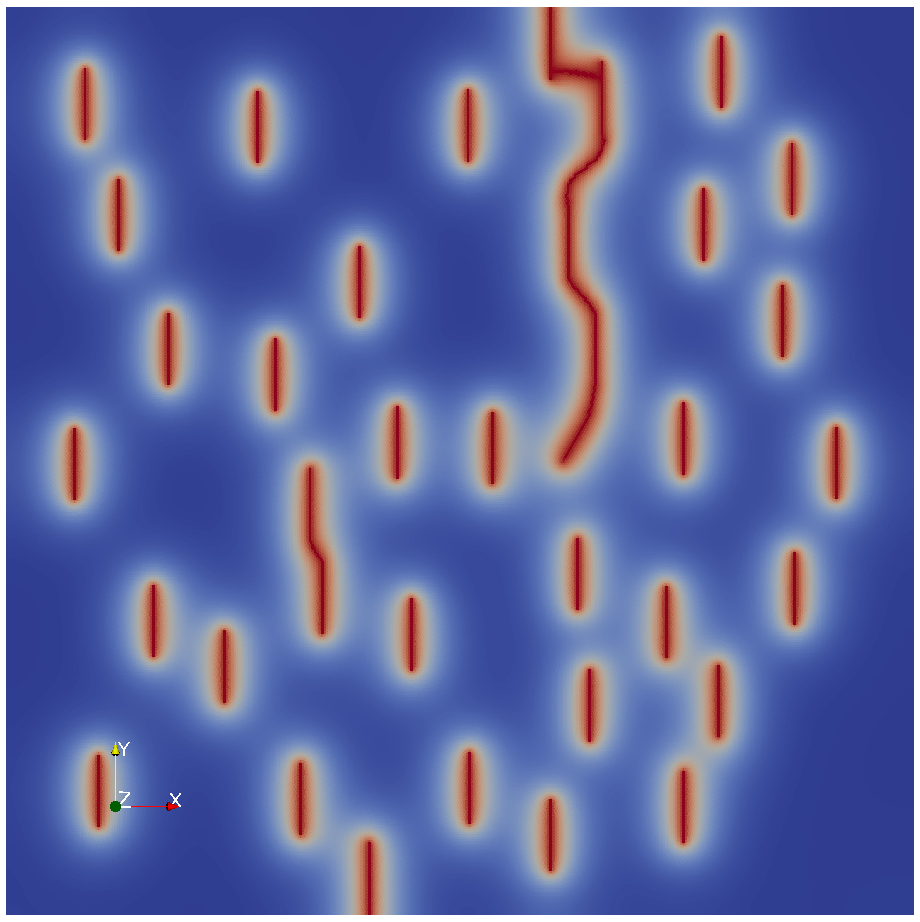}
		\caption{}
	\end{subfigure}
	\begin{subfigure}{0.3\textwidth}
		\includegraphics[width=\textwidth]{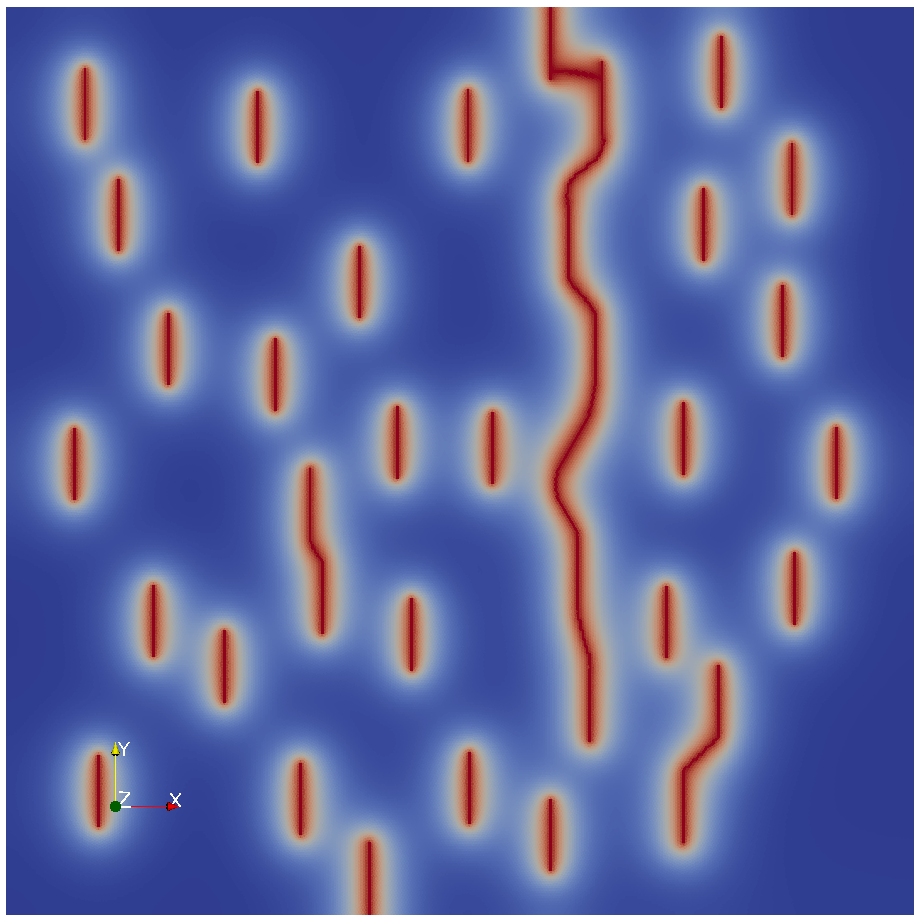}
		\caption{}
	\end{subfigure}
	\begin{subfigure}{0.3\textwidth}
		\includegraphics[width=\textwidth]{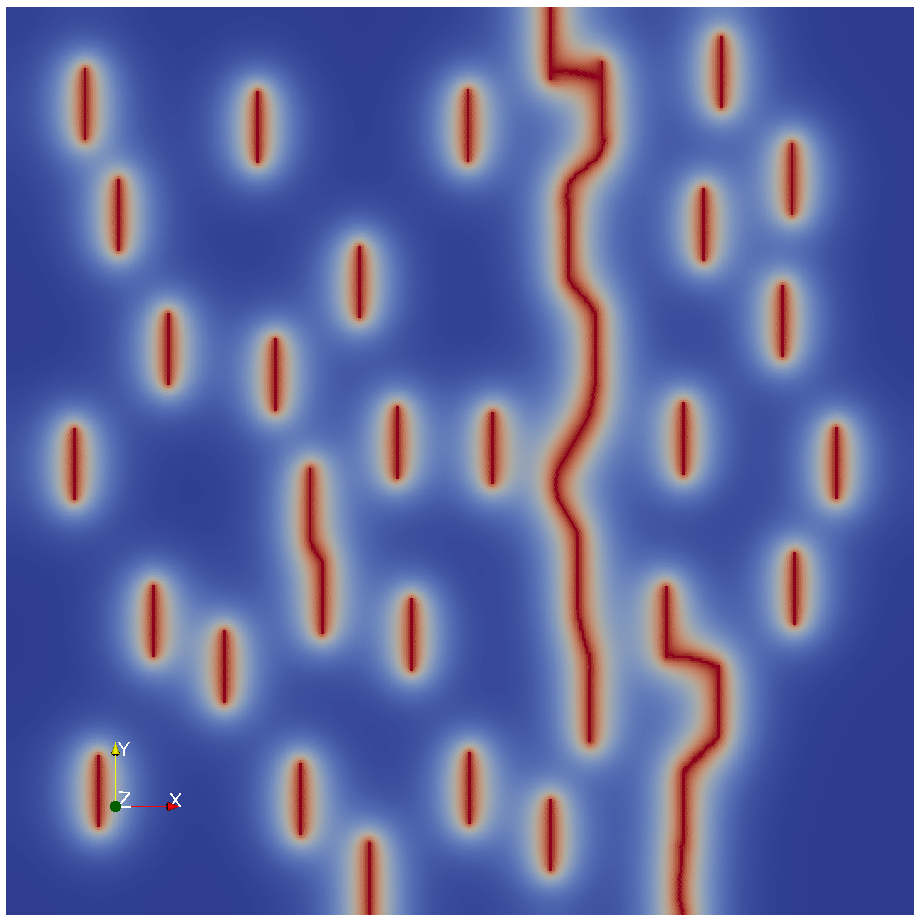}
		\caption{}
	\end{subfigure}
	\caption{Crack evolution for the specimen of Example \ref{sec:manyCracks}: (a) $u_x = 0.00613$, (b) $u_x = 0.00614$, (c) $u_x = 0.00625$, (d) $u_x = 0.00631$, (e) $u_x = 0.00632$, (f) $u_x = 0.00683$.}
	\label{fig:manyCrackEvolution}
\end{figure}

For the final example, we simulate fracture evolution under tensile loading in a square specimen having several preexisting cracks. The outer dimensions of the said specimen are 1 mm $\times$ 1 mm, and it is assumed to be made of the same material as those from the previous two examples. The initial cracks are scattered throughout the specimen; all are 0.08 mm long and oriented in the vertical direction with one crack intersecting the top boundary of the specimen and another the bottom as shown in Figure \ref{fig:manyCracksGeometry}.
For this problem, the phase-field regularization parameter is set to $\ell = 0.025$ and the domain is discretized uniformly into 245,326 elements having characteristic size $h = \ell/8$. The preexisting cracks are modeled through the history field that is initialized following the procedure outlined in Section \ref{sec:histFieldInitialization}, with the peak phase-field value set to 0.999. Furthermore this is done in a way such that the fully damaged region for each crack has a thickness of just one cell, so that the phase-field profile perpendicular to a crack is similar to that shown in Figure \ref{fig:comparison_1D}. The top and bottom of the specimen are assumed to be traction free, and its entire left boundary is fixed (i.e. $u_x = u_y = 0$). On the other hand right boundary is fixed in the vertical direction ($u_y = 0$) but subject to a horizontal displacement of $u_x = 0.01$ mm. This last BC is applied incrementally, with increment sizes having been manually adjusted in order to better capture the resulting fracture evolution. In the final run shown in the paper, said loading increments are as follows: $\Delta u_x = 1.0\times 10^{-4}$ from $u_x = 0$ to $u_x = 0.0057$, followed by $\Delta u_x = 1.0 \times 10^{-5}$ up to $u_x = 0.0069$, and thereafter $\Delta u_x = 1.0\times 10^{-4}$ until the target displacement of $u_x = 0.01$ is reached. For this problem we have opted to use the anisotropic model of Amor et al.\ \cite{Amor2009} in tandem with the single-parameter degradation function
\begin{linenomath}
\begin{equation}
	g_s \left( \phi; n \right) = \frac{1 - e^{-k \left( n \right) \left( 1 - \phi \right)^n}}{1 - e^{-k \left( n \right)}}
\end{equation}
\end{linenomath}
introduced in \cite{Sargado2018} to minimize spurious damage evolution in the pre-fracture elastic response that is known to occur when the traditional quadratic degradation function is used. In particular, we set $n = 2.0$ in the above expression to confine highly damaged regions to the near-crack vicinity.\footnote{For details on how to calculate $k \left( n \right)$, we direct interested readers to \cite{Sargado2018}. We note also that the complete model given therein has an additional term that is weighted by a parameter $w$; the latter has been set to zero for the current application.} The resulting load-displacement curve is shown in Figure \ref{fig:loadDisp_manyCracks}, where for reference we also plot the ``elastic'' solution utilizing a discontinuous phase-field, i.e. $\phi = 0.999$ in the critical region containing the initial cracks and 0 elsewhere.
We note that while there is indeed no appreciable gradual loss of stiffness in the mechanical response prior to fracture (i.e. the pre-fracture elastic behavior is linear), there is quite a significant discrepancy between the initial slope of the load-displacement curve for the evolving crack problem and that of the elastic solution. This is expected and can be attributed to the modified stress-strain behavior near crack tips as discussed in \cite{Sicsic2013}. Said discrepancy becomes worse with the number of initial cracks (since this means there are more crack tips) but can be made smaller by further reducing the value of $\ell$. However doing so can easily lead to prohibitive computational costs unless an adaptive remeshing strategy is used. Alternative approaches may also be possible such as choosing $\ell/h$ nearer to 1 and then compensating for the resulting error in the numerically predicted $\mathcal{G}_c$, but such approaches have not yet been rigorously investigated in the literature and are likewise beyond the scope of the current work. Phase-field profiles corresponding to specific points of the load-displacement curve are shown in Figure \ref{fig:manyCrackEvolution}.
\section{Concluding remarks}
In this paper we have introduced an unequal-order discretization scheme for brittle fracture phase-field models that combines $P_1$ finite element and cell-centered ($P_0$) finite volume approaches. This was motivated by the inference that mesh size restrictions with respect to the phase-field length scale may be much more stringent than the requirements to obtain sufficiently accurate solution of the stresses due to the fact that introduction of damage at the crack tips eliminates stress singularities. Application of the said formulation to two popular benchmark problems in the literature have shown that the computational cost of obtaining mesh-insensitive results are significantly reduced in comparison with an equal-order ($P_1$) finite element discretization of the mechanics and phase-field equations. This is rather counter-intuitive looking at the degrees of polynomial approximation, but may be explained by the fact that the manner in which flux/gradients are calculated across cell faces in the FV framework allows for the implicit occurrence of a cusp inside the control volume and can thus capture the phase-field profile more efficiently. As virtually all numerical implementations of fracture phase-field models in the literature utilize equal-order formulations, we consider this work as a proof of concept regarding the potential of using unequal-order and combined discretization approaches, for instance high order finite elements in tandem with discontinuous Galerkin approximations. Likewise it would be interesting to investigate the performance of the proposed formulation in 3D (i.e., $P_1$-$P_0$ based on linear tetrahedra).

% =====================================================
\section*{Acknowledgments}
This work was funded by the Research Council of Norway through grant no.\ 228832/E20 and Equinor ASA through the Akademia agreement.

% =====================================================
\appendix
\section{Implementation of strain decomposition for anisotropic bulk degradation} \label{sec:amorDetails}
In the code utilized for the current work, the isotropic and anisotropic formulations given in \eqref{eq:bourdinModel} and \eqref{eq:amorModel} are implemented as different constitutive models in conjunction with a single numerical framework that corresponds to the discretized field equations. These in turn utilize predefined modulus tensors for linear elasticity corresponding to plane strain, written in Voigt form as
\begin{linenomath}
\begin{equation}
	\mathbb{C}^e = \left[ \begin{array}{cccc}
		\lambda + 2\mu & \lambda & \lambda & 0 \\
		\lambda & \lambda + 2\mu & \lambda & 0 \\
		\lambda & \lambda & \lambda + 2\mu & 0 \\
		0 & 0 & 0 & \mu
	\end{array} \right],
\end{equation}
\end{linenomath}
where $\lambda$ and $\mu$ are the Lam\'e coefficients. As shown above, the matrix dimensions of $\mathbb{C}^e$ are $4\times 4$ since for applications involving plasticity, the $z$-component of plastic strain may be nonzero even if $\epsilon_{zz} = 0$. On the other hand in the bulk degradation model of Amor et al.~\cite{Amor2009}, the strain tensor is split into volumetric/spherical and deviatoric components, with the former used as a reference quantity for modeling mode I fracture and unilateral contact, and the latter for cracking in mode II. Thus in plane strain problems, one can argue that the deviatoric component of strain in the $z$-direction should not play a role in mode I/II fracture evolution as it instead implies cracking in mode III. Consequently, we perform the strain decomposition considering only two dimensions. Writing the strain in Voigt form as $\bm{\upepsilon} = \left\{ \epsilon_{xx}, \epsilon_{yy}, 0, \gamma_{xy} \right\}^\trp$, we calculate volumetric and deviatoric components as
\begin{linenomath}
\begin{equation}
\begin{split}
	\bm{\upepsilon}_\text{vol} &= \mathbb{P} \bm{\upepsilon} \\
	\bm{\upepsilon}_\text{dev} &= \left( \mathbb{I} - \mathbb{P} \right) \bm{\upepsilon}
\end{split}
\end{equation}
\end{linenomath}
where
\begin{linenomath}
\begin{equation}
	\mathbb{P} = \frac{1}{2} \left[ \begin{array}{cccc} 
		1 & 1 & 0 & 0 \\ 1 & 1 & 0 & 0 \\ 0 & 0 & 0 & 0 \\ 0 & 0 & 0 & 0
	\end{array} \right] \text{ and }
	\mathbb{I} = \left[ \begin{array}{cccc} 
	1 & 0 & 0 & 0 \\ 0 & 1 & 0 & 0 \\ 0 & 0 & 1 & 0 \\ 0 & 0 & 0 & 1
	\end{array} \right].
\end{equation}
\end{linenomath}
The above formula yields $\epsilon_{\text{vol},zz} = \epsilon_{\text{dev},zz} = 0$, hence $\bm{\upepsilon}_\text{vol}$ is no longer a spherical tensor. Denoting by $H \left( \bullet \right)$ the Heaviside step function, the tangent modulus accounting for damage is given by
\begin{linenomath}
\begin{equation}
	\mathbb{C} \left( \bm{\upepsilon}, \phi \right) = g \left( \phi \right) \mathbb{C}^e  \left\{ \left[ H \left( \trace \bm{\upepsilon} \right) - 1 \right] \mathbb{P}  + \mathbb{I} \right\} + H \left( -\trace \bm{\upepsilon} \right) \mathbb{C}^e \mathbb{P}
\end{equation}
\end{linenomath}
which will not be symmetric due to the form of $\mathbb{P}$. Nevertheless, the product $\mathbf{B}^\trp \mathbb{C} \left( \bm{\upepsilon}, \phi \right) \mathbf{B}$ does yield a symmetric matrix as the third row of $\mathbf{B}$ is composed entirely of zeros.

\section*{References}
\bibliographystyle{unsrt}
\bibliography{references}
\end{document}